\documentclass[english,a4paper,eqno,11pt]{article}

\usepackage{amssymb}
\usepackage{latexsym}
\usepackage[english]{babel}
\font\tenbb=msbm10 at 12pt

\def\rR{\hbox{\tenbb R}}
\def\nN{\hbox{\tenbb N}}

\def\cal{\mathcal}
\font\titre cmbx10 at 18 pt

\def\gesp{\vskip1cm}
\def\esp{\vskip .6cm}
\def\pesp{\vskip .3cm}

\def\ni{\noindent}
\def\di{\displaystyle}

\def\Box{$\sqcap\hskip-.262truecm\sqcup$}
\newtheorem{thm}{Theorem}[section]
\newtheorem{defn}{Definition}[section]
\newtheorem{lem}{Lemma}[section]
\newtheorem{prop}{Proposition}[section]
\newtheorem{prp}{Property}[section]
\newtheorem{cor}{Corollary}[section]
\newtheorem{Exp}{Example}[section]
\newtheorem{rem}{Remark}[section]

\usepackage{graphicx,epsfig,epic} % standard LaTeX graphics tool
\psfigdriver{dvips}

\usepackage{graphicx}
\usepackage{amscd}
\usepackage{latexsym,amssymb,euscript}
\usepackage{epic,eepic}
 \usepackage{pifont}
\begin{document}

\centerline{\titre {Fractal Topology Foundations }}

\esp
\centerline{ H\'EL\`ENE PORCHON}
\pesp
\centerline{\small UFR 929 Math\'ematiques }
\centerline{\small Universit\' e Pierre et Marie Curie}
\centerline{\small  4 Place Jussieu, 75252 Paris Cedex 05, France}
\centerline{helene.porchon@etu.upmc.fr}

\gesp

{\small Abstract. In this paper, we introduce the foundation of a fractal topological space constructed via
a family of nested topological spaces endowed with subspace topologies, where the number of topological
spaces involved in this family is related to the appearance of new structures on it. The greater the number of topological spaces we use,
the stronger the subspace topologies we obtain. The fractal manifold model is brought up as an illustration of space that is locally homeomorphic to the fractal topological space.
 }\pesp

{\small Keywords}: {\small Differentiable Manifold; Topological Space; Fractal.}

{\small MSC (2010)}: {\small 54A05 - 54A10 - 54D80 - 54F65 - 54H20}

\section{Introduction}

The notion of "fractal topology" is usually used in various domains as for example cosmic electrodynamic \cite{ZM} or complex network from different fields
(biology, technology, sociology) \cite{SH}.
In a mathematical point of view, attempts have been done to relate topology and fractal geometry \cite{Ed}, or to describe the topology of fractal sets \cite{MAZ}.
However there does not exist any mathematical foundation and formulation for a fractal topology.

The concept of fractal topology presented in this paper is a new concept derived from the fractal manifold model  \cite{BF1}.
Indeed a fractal manifold is typically the kind of space that is naturally endowed with a fractal
topology.

This paper presents in a preliminary part (section 2) an introduction to basic notions related to the fractal manifold and its properties \cite{BF1},\cite{BF3}.
The main results are given in a second part (section 3): we first give a definition of a fractal family of topological spaces, and of a fractal topology. Then we
study the fractal manifold to determine its topology.
\pesp

\section{Preliminary}

We introduce in this part basic notions about $\delta_0$-manifold, diagonal topology
and fractal manifold that can be found in \cite{BF1} with deeper details related to the construction.\pesp

\subsection{${\delta_0}$-Manifold}

Let $f_i$, for $i=1,2,3,$ be three continuous and nowhere
differentiable functions,
defined on the interval $[a,b]\subset\rR$, with $a<b$. For $i=1,2,3$, the associated graph of $f_i$ is given by
$\Gamma_{i,0}([a,b])=\Big\lbrace(x,y)\in\rR^2/\ y=f_i(x),
\ x\in[a,b]\Big\rbrace.$
For $i=1,2,3,$ let us consider the function $f_i(x,y)={1\over 2y}\int^{x+y}_{x-y}f_i(t)dt,$ we call forward (respectively backward) mean
function of $f_i$ the function given by:
\begin{equation}\label{E0}
f_i(x+\sigma_0{\delta_0\over2},{\delta_0\over2})=\di{\sigma_0\over\delta_0}\int_x^{x+\sigma_0\delta_0}f_i(t)dt
\quad \hbox{for}\ \sigma_0=+\ (\hbox{respectively}\ \sigma_0=-),\hbox{and}\  \delta_0 \in \rR,
\end{equation}
and we denote by
$\Gamma_{i,\delta_0}^{\sigma_0}$ its associated graph.
\pesp
We define the translation
$T_{\delta_0} :
\prod_{i=1}^{3}\Gamma_{i{\delta_0}}^{+}\times \{{\delta_0}\}
\longrightarrow \prod_{i=1}^{3}\Gamma_{i{\delta_0}}^{-}\times
\{{\delta_0}\}$ by
$$T_{\delta_0}
\Big((a_1,b_1),(a_2,b_2),(a_3,b_3)\Big)=\Big((a_1+{\delta_0},b_1),(a_2+{\delta_0},b_2),
(a_3+{\delta_0},b_3)\Big),$$
where $(a_i,b_i)\in
\Gamma_{i{\delta_0}}^{+}$, that is to say\quad  $b_i=\di
f_i(a_i+{{\delta_0} \over 2}, {{\delta_0} \over 2})={1\over
{\delta_0}}\di \int _{a_i}^{a_i+{\delta_0}}f_i(t)dt$
for   $i=1,2,3$.
\pesp
Let us consider $\varepsilon_n\in]0,1[$ for all $n\geq0$ such that $\varepsilon_0>\varepsilon_1>\ldots>\varepsilon_n$, $\forall n>0$.
We denote ${\cal R}_{0}=]0,\varepsilon_0]$, and  ${\cal R}_{n}=[0,\varepsilon_n]$ for all $n>0$.
In the following we consider the nested real numbers $\delta_0>\delta_1>\ldots>\delta_n$, $\forall n>0$, and such that $\delta_n\in{\cal R}_{n}$ for all $n\geq0$.
\pesp

\begin{defn}\label{Def0}
For $\delta_0\in{\cal R}_{0}$, let ${\cal M}_{\delta_0}$ be an
Hausdorff topological space. We say that ${\cal M}_{\delta_0}$ is an
${\delta_0}$-manifold if for every point $x \in M_{\delta_0}$,
there exist a neighborhood $\Omega _{{\delta_0}}$ of $x$ in
$M_{\delta_0}$, a map $\varphi _{\delta_0}$, and two open sets
$V^{+} _{{\delta_0}}$ of
$\prod_{i=1}^{3}\Gamma_{i{\delta_0}}^{+}\times \{{\delta_0}\} $
and $V^{-} _{{\delta_0}}$ of
$\prod_{i=1}^{3}\Gamma_{i{\delta_0}}^{-}\times \{{\delta_0}\} $
such that $\varphi _{\delta_0} : \Omega _{{\delta_0}}
\longrightarrow V^{+} _{{\delta_0}}$, and $T_{\delta_0} \circ
\varphi _{\delta_0}: \Omega _{{\delta_0}} \longrightarrow V^{-}
_{{\delta_0}} $ are two homeomorphisms.
\end{defn}\pesp

\begin{rem}
A ${\delta_0}$-manifold ${\cal M}_{{\delta_0}}$ is locally seen through a triplet $(\Omega _{{\delta_0}},\varphi _{\delta_0},T_{\delta_0}\circ
\varphi _{\delta_0})$ as illustrated in Diagram A, that is to say a point of ${\cal M}_{{\delta_0}}$ is represented in the local chart by two points $x^+$ and $x^-$
respectively in $\prod_{i=1}^{3}\Gamma_{i{\delta_0}}^{+}\times \{{\delta_0}\} $ and $\prod_{i=1}^{3}\Gamma_{i{\delta_0}}^{-}\times \{{\delta_0}\}$, and a
neighborhood in ${\cal M}_{{\delta_0}}$ is seen in $\prod_{i=1}^{3}\Gamma_{i{\delta_0}}^{+}\times \{{\delta_0}\} $ and $\prod_{i=1}^{3}\Gamma_{i{\delta_0}}^{-}\times
\{{\delta_0}\}$ respectively as two neighborhoods $V^{+} _{{\delta_0}}$ and $V^{-} _{{\delta_0}}$, where $V^{-} _{{\delta_0}}$ is obtained from $V^{+} _{{\delta_0}}$
by the translation $T_{\delta_0} $.
\unitlength=1.2cm
\begin{picture}(11,4.5)

%%%%%++++

\put(6,4){$\prod_{i=1}^{3} \Gamma_{i{\delta_0}}^{+}\times \{{\delta_0}\}$}

%%%%%-----

\put(6,1.5){$\prod_{i=1}^{3} \Gamma_{i{\delta_0}}^{-}\times \{{\delta_0}\}$}

% les vecteurs
\put(4,3){\vector(2,1){1.6}}

\put(6.1,3.5){\vector(0,-3){1.5}}

\put(4,2.5){\vector(2,-1){1.5}}

%les fonctions
\put(4.5,3.7){$\varphi_{\delta_0}$}

\put(4.1,1.6){$T_{\delta_0}\circ\varphi_{\delta_0}$}

\put(6.3,2.7){$T_{\delta_0}$}

%la variete
\put(3,2.7){${\cal M}_{\delta_0}$}
%%%%%%%%%%%%%%%%%%%%%%%
\put(1.5,0.3){\footnotesize { Diagram A: $\delta_0$-manifold locally defined by a double homeomorphism.}}
\thicklines
\end{picture}
\end{rem}
\pesp

\subsection{Diagonal Topology}

In the purpose to define a fractal manifold, we will use the notion of "diagonal topology" introduced in \cite{BF1}.
Let us consider in general a set $E=\cup_{\varepsilon\in I}E_\varepsilon $ union of topological spaces all disjoint or all the same, where ${\cal T}_\varepsilon$
is the topology on $E_\varepsilon$ for all $\varepsilon\in I$, and where $I$ is a bounded interval of $\rR$.

\begin{prp}
If $A=\cup_{\varepsilon\in I}A_\varepsilon $ and $B=\cup_{\varepsilon\in I}B_\varepsilon $ are two subsets of $E=\cup_{\varepsilon\in I}E_\varepsilon $
such that $A_\varepsilon $ and $B_\varepsilon $ are subsets of $E_\varepsilon $ for all $\varepsilon\in I$, then
$A{\cap}B=\cup_{\varepsilon\in I}(A_\varepsilon\cap B_\varepsilon)$ and $A{\cup}B=\cup_{\varepsilon\in I}(A_\varepsilon\cup B_\varepsilon)$, where
$A_\varepsilon\cap B_\varepsilon\subset E_\varepsilon$ and $A_\varepsilon\cup B_\varepsilon\subset E_\varepsilon$ for all $\varepsilon\in I$.
\end{prp}
\pesp

We can now define  a diagonal topology on $E=\cup_{\varepsilon\in I}E_\varepsilon$:

\begin{defn}\label{DefTop}
We call diagonal topology on $E=\cup_{\varepsilon\in I}E_\varepsilon $ union of topological spaces all disjoint or all the same, the topology
${\cal T}_d$ defined by
$${\cal T}_d= \Big\{ \Omega=\cup_{\varepsilon\in I}\Omega_\varepsilon \quad / \ \Omega_\varepsilon\in
{\cal T}_\varepsilon\quad \forall \varepsilon\in I\ \Big\}$$
where ${\cal T}_\varepsilon$ is the topology on $E_\varepsilon$ for all $\varepsilon\in I$. The topological space
$(E, {\cal T}_d)$ is called diagonal topological space.
\end{defn}\pesp

\begin{rem}
The diagonal topology ${\cal T}_d$ is a topology on $E$ since it satisfies the following axioms:

i) $\emptyset \in {\cal T}_d$ and $E\in {\cal T}_d$

ii) if $\ \Omega_1 \in {\cal T}_d\ $ and $\ \Omega_2 \in {\cal T}_d$,\  then \ \ $\Omega_1\ {\cap}\ \Omega_2\in {\cal T}_d$

iii) any union of elements of ${\cal T}_d$ is an element of ${\cal T}_d$.
\end{rem}
\pesp

In the following, we will use the diagonal topology ${\cal T}_d$ on any union of topological spaces all disjoint or all the same. Therefore we need some
specific definitions related to the topology ${\cal T}_d$.

\begin{defn}\label{Def1}
We call object of $(E,{\cal T}_d)$ a set $X=\cup_{\varepsilon\in I}\{x_\varepsilon\}$, where $x_\varepsilon\in E_\varepsilon$ for all $\varepsilon\in I$.
\end{defn}\pesp

Therefore an object of $(E,{\cal T}_d)$ is a family of points that has a representative element in each $E_\varepsilon$, $\varepsilon\in I$.

\begin{defn}
Let us consider an object $X=\cup_{\varepsilon\in I}\{x_\varepsilon\}$ of $(E,{\cal T}_d)$. A subset $\omega$ of $E$ is called a diagonal neighborhood of
$X$ if there exists $\Omega=\cup_{\varepsilon\in I}\Omega_\varepsilon \in {\cal T}_d$ such that $\Omega\subset\omega$ and  $\Omega_\varepsilon $ is a
neighborhood of $x_\varepsilon$ for all $\varepsilon\in I$.
\end{defn}\pesp

\begin{defn}\label{Def8}
We say that $E$ admits an internal structure $\tilde{X}$ on an element $P\in E$ if there exists a $C^0$-parametric path
\begin{equation}
\begin{array}{lll}
\tilde{X}: I & \longrightarrow & \bigcup_{\varepsilon\in {I}} E_{\varepsilon} \\
\quad\ \  {\varepsilon} & \longmapsto & \tilde{X}(\varepsilon)\in E_{\varepsilon} ,
\end{array}
\end{equation}
such that for all ${\varepsilon}\in I$, $Range (\tilde{X}) \cap
E_{\varepsilon}=\Big\{\tilde{X}({\varepsilon})\Big\}\ $, and there exists ${\varepsilon}'\in I$ such that
$P=\tilde{X}(\varepsilon')\in E_{\varepsilon'}$. We call $P$ a point of $E$.
\end{defn}\pesp

\begin{rem}
The set $Range (\tilde{X})=\cup_{\varepsilon\in I}\{\tilde{X}(\varepsilon)\}$ is an object of $E$ since $\tilde{X}(\varepsilon)\in E_\varepsilon\
\forall \varepsilon\in I $.
\end{rem}\pesp

We can define on the set of internal structures of $E$ an equivalence relation which allows to talk about uniqueness of the internal structure.

\begin{defn}\label{Def2}
Let $\tilde{X}$ and $\tilde{Y}$ be two internal structures of $E$. We say that $\tilde{X}\sim \tilde{Y}$ if and only if

i) $\exists\ \varepsilon'\in I$ such that $\tilde{X}(\varepsilon')=\tilde{Y}(\varepsilon')$.

ii) $\exists\ \theta:I \rightarrow I $ diffeomorphism such that $\tilde{X}=\tilde{Y}\circ \theta$.
\end{defn}\pesp

We send the reader to \cite{BF1} for the proof of the following proposition:

\begin{prop}
$\tilde{X} \sim \tilde{Y} \Leftrightarrow \tilde{X}=\tilde{Y}$.
\end{prop}\pesp

\subsection{Fractal-Manifold}

We consider now an union ${\cal M}=\bigcup_{{\delta_0}\in {\cal R }_0}{\cal M}_{{\delta_0}}$ of $\delta_0$-manifolds all disjoint or all the same, where the variable
$\delta_0$ varies in ${\cal R }_0$.\pesp

\begin{defn}\label{Def7}
A fractal manifold is an union of Hausdorff topological spaces all disjoint or all the same ${\cal M}=\bigcup_{{\delta_0}\in {\cal R }_0}{\cal M}_{{\delta_0}}$,
which satisfies the following properties:
$\forall{\delta_0} \in {\cal R}_0$,
${\cal M}_{\delta_0}$ is a ${\delta_0}$-manifold,
 and  $\forall P\in {\cal M}$, ${\cal M}$ admits an internal structure $\tilde{X}$ on $P$ such
that there exist a neighborhood\quad
 $\Omega(Rg(\tilde{X}))=\cup_{{\delta_0}\in {\cal R}_0}\Omega_{\delta_0}$,
with $\Omega_{\delta_0}$ a neighborhood of $\tilde{X}({\delta_0})$ in
$M_{\delta_0}$, two open sets $V^+=\cup_{{\delta_0}\in {\cal
R}_0} V_{\delta_0}^+$ and $V^-=\cup_{{\delta_0}\in {\cal R}_0}
V_{\delta_0}^-$, where $V_{\delta_0}^\sigma$ is an open set in
$\Pi_{i=1}^3\Gamma_{i{\delta_0}}^\sigma\times\{{\delta_0}\}$ for
$\sigma=\pm$, and there exist two families of maps
$(\varphi_{\delta_0})_{{\delta_0}\in {\cal R}_0}$ and
$(T_{\delta_0}\circ\varphi_{\delta_0})_{{\delta_0}\in {\cal
R}_0}$ such that
$\varphi_{\delta_0}:\Omega_{\delta_0}\longrightarrow
V_{\delta_0}^+ $ and
$T_{\delta_0}\circ\varphi_{\delta_0}:\Omega_{\delta_0}\longrightarrow
V_{\delta_0}^-$ are homeomorphisms for all ${{\delta_0}\in {\cal
R}_0}$.
\end{defn}\pesp

\begin{defn}\label{Def6}
i) A local chart on the fractal manifold ${\cal M}$ is a triplet $(\Omega,
\varphi, T \circ \varphi)$, where $\Omega=\bigcup _{{\delta_0}
\in {\cal R}_0}\Omega_{{\delta_0}}$ is an open set of ${\cal M}$,
$\varphi= (\varphi_{\delta_0})_{{\delta_0}\in {\cal R}_0}$ is a family of homeomorphisms $\varphi_{\delta_0}$
from $\Omega_{{\delta_0}}$ to an open set $V^+_{\delta_0}$ of\quad
$\prod_{i=1}^{3}\Gamma_{i{\delta_0}}^{+}\times \{{\delta_0}\}$,
and $T\circ \varphi=(T_{\delta_0}\circ \varphi_{\delta_0})_{{\delta_0}\in {\cal R}_0}$ is a family of homeomorphisms
$T_{\delta_0}\circ \varphi_{\delta_0}$ from
$\Omega_{{\delta_0}}$ to an open set $V^-_{\delta_0}$ of\quad
$\prod_{i=1}^{3}\Gamma_{i{\delta_0}}^{-}\times \{{\delta_0}\} $
for all ${\delta_0}\in{\cal R}_0$.

ii) A collection $\ (\Omega_{i},
\varphi_{i}, (T\circ \varphi)_{i})_{i\in J}\ $ of local charts on
the fractal manifold ${\cal M}$ such that $\ \Omega_{i}=\bigcup _{{\delta_0}
\in {\cal R}_0}\Omega_{i,{\delta_0}}$, $\ \varphi_{i}=(\varphi_{i,\delta_0})_{{\delta_0}\in {\cal R}_0}$, $\ (T\circ\varphi)_{i}=(T_{i,{\delta_0}}\circ
\varphi_{i,\delta_0})_{{\delta_0}\in {\cal R}_0}\ $
and such that $\ \cup _{i\in
J}\Omega_{i,{\delta_0}}=M_{\delta_0}$ for all ${\delta_0}\in{\cal R}_0$, and $\cup _{i\in J}\Omega _i=M$
is called an atlas.

iii) The coordinates of an object $Rg(\tilde{X})\subset \Omega$ related to the local
chart $(\Omega, \varphi, T \circ \varphi )$ are the coordinates of
the object $\ \varphi(Rg(\tilde{X}))\ $ in $\ \bigcup _{{\delta_0} \in {\cal
R}_0}\prod_{i=1}^{3}\Gamma_{i{\delta_0}}^{+}\times
\{{\delta_0}\}\ $, and of the object\quad $(T \circ \varphi)(Rg(\tilde{X}))$ in
$\bigcup _{{\delta_0} \in {\cal
R}_0}\prod_{i=1}^{3}\Gamma_{i{\delta_0}}^{-}\times
\{{\delta_0}\} $.
\end{defn}\pesp

\begin{rem}
1) An illustration of the fractal manifold ${\cal M}$ is given by the Diagram B:

\unitlength=1.2cm
\begin{picture}(11,4.7)

%%%%%++++

\put(6.5,3.9){$\bigcup_{{\delta_0}\in {\cal R}_0}\prod_{i=1}^{3} \Gamma_{i{\delta_0}}^{+}\times \{{\delta_0}\}$}

%%%%%-----

\put(6.5,1.5){$\bigcup_{{\delta_0}\in {\cal R}_0}\prod_{i=1}^{3} \Gamma_{i{\delta_0}}^{-}\times \{{\delta_0}\}$}

% les vecteurs
\put(3.5,3){\vector(3,1){2.5}}

\put(7.6,3.5){\vector(0,-3){1.5}}

\put(3.5,2.5){\vector(3,-1){2.5}}

%les fonctions
\put(2.2,3.6){\small $\varphi_1=(\varphi _{\delta_0})_{\delta_0 \in {\cal R}_0}$}

\put(1.6,1.7){\small $T_1\circ\varphi_1=(T_{\delta_0}\circ \varphi _{\delta_0})_{\delta_0 \in {\cal R}_0}$}

\put(7.8,2.7){\small $T_1=(T_{\delta_0})_{\delta_0 \in {\cal R}_0}$}

%la variete
\put(0.7,2.7){${\cal M}=\bigcup_{{\delta_0}\in {\cal R}_0}{\cal M}_{\delta_0}$}

\put(0.2,0.5){ \footnotesize { Diagram B: Fractal manifold ${\cal M}$ locally defined by a double family of homeomorphisms.}}
\thicklines
\end{picture}\pesp

2) We have on the set of internal structures of ${\cal M}$ the equivalence relation $\sim$ defined in Definition \ref{Def2}. Moreover ${\cal M}$ admits a unique internal
structure at each point by the definition of fractal manifold, then ${\cal M}$ can be seen as the set of all equivalence classes for the equivalence relation $\sim$,
and points of ${\cal M}$ can be assimilated to objects of ${\cal M}$.

3) By Definition \ref{Def1}, if
$\tilde{X}:\ {\cal R}_f\subset\rR\ \longrightarrow\ \cup_{{\delta_0}\in {\cal R}_0}
{\cal M}_{{\delta_0}}$ is an internal structure on ${\cal M}=\bigcup_{{\delta_0}\in {\cal R }_0}{\cal M}_{\delta_0}$, then the set $Rg(\tilde{X})=\cup_{{\delta_0}\in {\cal R}_0}
\{\tilde{X}(\delta_0)\}$ is an object of ${\cal M}$.
\end{rem}\pesp

For sake of simplicity, let us introduce the following notations.\pesp

\ni{\bf Notations:}\pesp

1) For all $ n\geq0$, for all $\delta_0\in {\cal R}_0,\ldots,\delta_n\in{\cal R}_n$ and $\sigma_0=\pm,\ldots,\sigma_n=\pm$, we denote by
$\di N_{\delta_0...\delta_{n}}^{\sigma_0...\sigma_{n}}$ the following set:

\begin{equation}\label{F10}
 N_{\delta_0...\delta_{n}}^{\sigma_0...\sigma_{n}}=\prod_{i=1}^3\Gamma_{i\delta_{0}...\delta_{n}}^{\sigma_0...\sigma_{n}}
\times\{\delta_{n}\}
\times...\times\{\delta_{0}\}
\end{equation}

\ni where
$\Gamma_{i\delta_{0}...\delta_{n}}^{\sigma_0...\sigma_{n}}$ represents the graph of the function:

\begin{equation}\label{mom}
F^{\sigma_0...\sigma_{n}}_{i\delta_0...\delta_{n}}(x)={\sigma_n...\sigma_0\over\delta_n...\delta_0}\int_x^{x+\sigma_n\delta_n}
\int_{t_{n-1}}^{t_{n-1}+\sigma_{n-1}\delta_{n-1}}\ldots\int_{t_{0}}^{t_0+\sigma_0\delta_0}f_i(t)dtdt_0\ldots dt_{n-1},
\end{equation}

\ni that is to say:\pesp

i) for $n=0$ we have 2 graphs $\Gamma_{i\delta_{0}}^{+}$ and $\Gamma_{i\delta_{0}}^{-}$.

ii) for $n=1$ we have 4 graphs $\Gamma_{i\delta_{0}\delta_{1}}^{++},\ \Gamma_{i\delta_{0}\delta_{1}}^{+-},\ \Gamma_{i\delta_{0}\delta_{1}}^{-+}$ and
$\Gamma_{i\delta_{0}\delta_{1}}^{--}$

iii) for $n=2$ we have 8 graphs $\Gamma_{i\delta_{0}\delta_{1}\delta_{2}}^{+++},\ \Gamma_{i\delta_{0}\delta_{1}\delta_{2}}^{++-},\ \Gamma_{i\delta_{0}\delta_{1}\delta_{2}}^{+-+},
\ \Gamma_{i\delta_{0}\delta_{1}\delta_{2}}^{+--},\ \Gamma_{i\delta_{0}\delta_{1}\delta_{2}}^{-++},\ \Gamma_{i\delta_{0}\delta_{1}\delta_{2}}^{-+-},
\ \Gamma_{i\delta_{0}\delta_{1}\delta_{2}}^{--+}$

 and $\Gamma_{i\delta_{0}\delta_{1}\delta_{2}}^{---}$.

iv) more generally for $n\geq0$, we have $2^{n+1}$ graphs.
\pesp

2) We denote the set $\di\bigcup_{\delta_0\in{\cal R}_0}\Big(\ldots\Big(\di\bigcup_{\delta_n\in{\cal R}_n}N_{\delta_0...\delta_{n}}^{\sigma_0...\sigma_{n}}\Big)\Big)$ by
$\di\bigcup_{\delta_0\ldots\delta_n}N_{\delta_0...\delta_{n}}^{\sigma_0...\sigma_{n}}$.
\esp

Using the previous notations, we introduce the following theorem (\cite{BF1}) that explains the internal chain react that gives
the fractal nature to the fractal manifold:\pesp

\begin{thm}\label{Th1}
If ${\cal M}$ is a fractal manifold,
then for all $n\geq0$, and for all $k\in [2^{n},2^{n+1}-1]\cap \nN$, there exist a family of local homeomorphisms
$\varphi_k$ and a family of translations $T_k$ such that for $\sigma_j=\pm$, $j=0,1,...,n$, one has the $2^{n}$ diagrams at the $step(n)$ given by Diagram C:
\gesp\esp

\unitlength=1.1cm
\begin{picture}(11,7)
\put(5.1,7.2){$ {\cal M}$}
%%%%%1111111
\put(3.5,5.7){\tiny$\bigcup_{\delta_0} N^+_{\delta_0}$}

\put(6,5.7){\tiny$\bigcup_{\delta_0} N^-_{\delta_0}$}

%%%%%%%2222222

\put(1.5,3.5){\tiny$\di\bigcup_{\delta_0\delta_1}
N^{++}_{\delta_0\delta_1}$}

\put(3.5,3.5){\tiny$\di\bigcup_{\delta_0\delta_1}
N^{+-}_{\delta_0\delta_1}$}

\put(5.9,3.5){\tiny$\di\bigcup_{\delta_0\delta_1}
N^{-+}_{\delta_0\delta_1}$}

\put(7.9,3.5){\tiny$\di\bigcup_{\delta_0\delta_1}
N^{--}_{\delta_0\delta_1}$}

%333333

\put(0.95,1){\tiny$\di\bigcup_{\delta_0\delta_1\delta_2}
N^{+++}_{\delta_0\delta_1\delta_2}$}

\put(2.5,-0.3){\tiny${\di\bigcup_{\delta_0\delta_1\delta_2}
N^{++-}_{\delta_0\delta_1\delta_2}}$}

\put(3,1){\tiny$\di\bigcup_{\delta_0\delta_1\delta_2}
N^{+-+}_{\delta_0\delta_1\delta_2}$}

\put(4.5,-0.3){\tiny${\di\bigcup_{\delta_0\delta_1\delta_2}
N^{+--}_{\delta_0\delta_1\delta_2}}$}

\put(5.3,1){\tiny$\di\bigcup_{\delta_0\delta_1\delta_2}
N^{-++}_{\delta_0\delta_1\delta_2}$}

\put(6.8,-0.3){\tiny${\di\bigcup_{\delta_0\delta_1\delta_2}
N^{-+-}_{\delta_0\delta_1\delta_2}}$}

\put(7.45,1){\tiny$\di\bigcup_{\delta_2\delta_1\delta_0}
N^{--+}_{\delta_0\delta_1\delta_2}$}

\put(9,-0.3){\tiny${\di\bigcup_{\delta_0\delta_1\delta_2}
N^{---}_{\delta_0\delta_1\delta_2}}$}

%les vecteurs obliques premiere generation
\put(5.2,7){\vector(-1,-1){1}}

\put(5.4,7){\vector(1,-1){1}}

% deuxieme Generation

\put(3.8,5.5){\vector(-1,-1){1.6}}

\put(4.1,5.5){\vector(0,-3){1.5}}

\put(6.4,5.5){\vector(0,-3){1.5}}

\put(6.8,5.5){\vector(1,-1){1.6}}

% troixieme Generation

\put(2.1,3.2){\vector(-1,-3){.6}}

\put(2.1,3.2){\vector(1,-3){1}}

%%%%%%%

\put(4.1,3.2){\vector(-1,-3){.6}}

\put(4.1,3.2){\vector(1,-3){1}}

%%%%%%%

\put(6.4,3.2){\vector(-1,-3){.6}}

\put(6.4,3.2){\vector(1,-3){1}}

%%%%%

\put(8.5,3.2){\vector(-1,-3){.6}}

\put(8.5,3.2){\vector(1,-3){1}}
%les vecteurs horizontale
\put(4.5,5.7){\vector(1,0){1.4}}

\put(2.5,4){\vector(1,0){1.5}}

\put(6.6,4){\vector(1,0){1.5}}

\put(1.4,.5){\vector(1,0){.9}}

\put(3.5,.5){\vector(1,0){.9}}

\put(5.7,.5){\vector(1,0){.9}}

\put(8,.5){\vector(1,0){.9}}
%les echelles
\put(10.74,6.5){\tiny$\delta_0$}

\put(10.74,4.5){\tiny$\delta_0,\delta_1$}

\put(10.74,1.7){\tiny$\delta_0,\delta_1,\delta_2$}

%\put(12.74,-1){\tiny$\delta_0,\delta_1,..,\delta_i,$}

%les fonctions

 \put(4.2,6.5){\tiny$\varphi_1$}

 \put(6.2,6.5){\tiny$T_1\circ\varphi_1$}
%%%%%%%%%%%%%%%%%%%%%%%%%%%%
 \put(2.5,4.8){\tiny$\varphi_2$}

 \put(4.3,4.8){\tiny$T_2\circ\varphi_2$}

 \put(6,4.8){\tiny$\varphi_3$}

 \put(7.9,4.8){\tiny$T_3\circ\varphi_3$}
%%%%%%%%%%%%%%%%%%%%%%%%%%%%%%%

\put(1.5,2.5){\tiny$\varphi_4$}

\put(2.5,2.1){\tiny$T_4\circ\varphi_4$}

\put(3.5,2.5){\tiny$\varphi_5$}

\put(4.5,2.1){\tiny$T_5\circ\varphi_5$}

\put(5.7,2.5){\tiny$\varphi_6$}

\put(6.8,2.1){\tiny$T_6\circ\varphi_6$}

\put(7.8,2.5){\tiny$\varphi_7$}

\put(8.9,2.1){\tiny$T_7\circ\varphi_7$}

\put(5.2,5.8){\tiny$T_{1}$}

\put(3.2,4.2){\tiny$T_{2}$}

\put(6.9,4.2){\tiny$T_{3}$}

\put(1.8,.2){\tiny$T_{4}$}

\put(3.8,.2){\tiny$T_{5}$}

\put(5.9,.2){\tiny$T_{6}$}

\put(8.3,.2){\tiny$T_{7}$}

%Titre
\put(1,-1.2){$\vdots\vdots$}

\put(2.6,-1.2){$\vdots\vdots$}

\put(3.3,-1.2){$\vdots\vdots$}

\put(5,-1.2){$\vdots\vdots$}

\put(5.7,-1.2){$\vdots\vdots$}

\put(7,-1.2){$\vdots\vdots$}

\put(7.7,-1.2){$\vdots\vdots$}

\put(9.5,-1.2){$\vdots\vdots$}
%line
\put(10.5,6){\line(0,1){1}}

\put(10.5,5.5){\bf $step(0)$}

\put(10.5,4){\line(0,1){1}}

\put(10.5,3.5){\bf $step(1)$}

\put(10.5,0.5){\line(0,1){2.5}}

\put(10.5,0){\bf $step(2)$}

\put(10.5,-1.2){ \vdots}

%%%%%%%%%%%%%%%%%%%%%%%
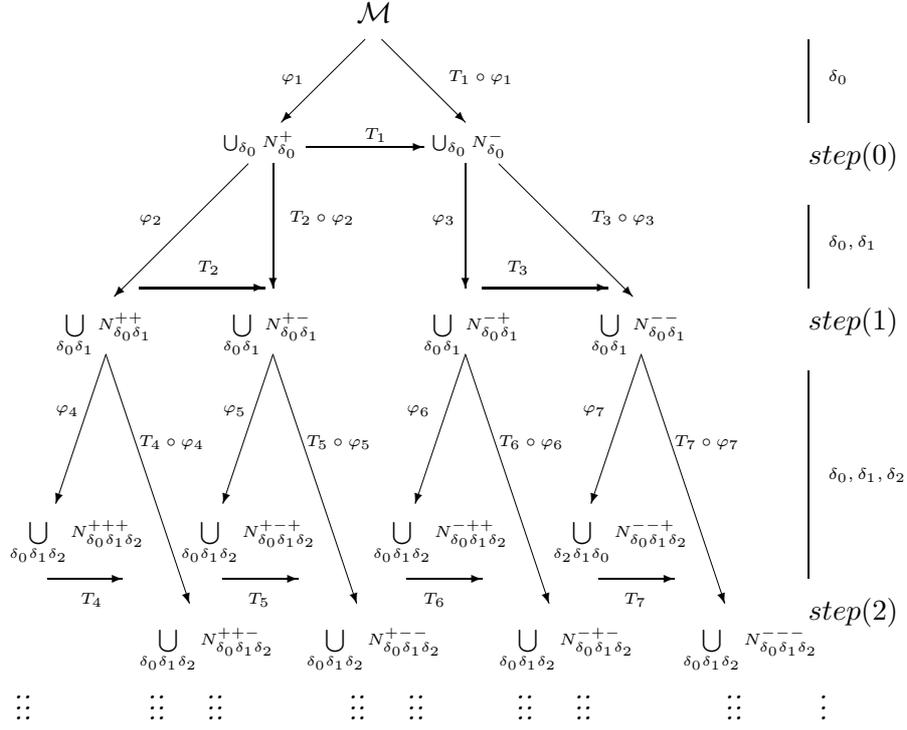
\put(2.5,-2){\footnotesize { Diagram C: Expanding diagram of the fractal manifold ${\cal M}$.}}

\thicklines
\end{picture}

\vskip3cm
\end{thm}
\pesp

\begin{rem}\label{R1}
The manifold ${\cal M}$ defined in Definition \ref{Def7} is called "Fractal Manifold" due to the theorem \ref{Th1} that rises the appearance of new structures at
each step (appearance of a new dimension $\delta_i$ at each step). The number of new dimensions increases as the number of steps increases.
\end{rem}\pesp

\section{Main Results}

The main objective is to prove that a fractal manifold has locally a fractal topology. Therefore we have to investigate the kind of structure necessary on a space
in order to precisely define the notion of fractal topology.

\subsection{Fractal Family of Topological Spaces and Fractal Topology}

The first task is to define the general concept of fractal topology. In this purpose we need to introduce the notion of fractal family of topological spaces.
\pesp

To compare diagonal topologies, it is convenient to introduce the following general definition:

\begin{defn}\label{ET}
Let $I$ be an interval of $\rR$, let $(E=\cup_{\varepsilon\in I}E_\varepsilon, {\cal T})$ and $(F=\cup_{\varepsilon\in I}F_\varepsilon, {\cal T}')$
be two diagonal topological spaces with diagonal topology respectively given by
$${\cal T}= \Big\{ \Omega=\cup_{\varepsilon\in I}\Omega_\varepsilon / \ \Omega_\varepsilon\in
{\cal T}_\varepsilon\quad \forall \varepsilon\in I\Big\} \quad\hbox{and} \quad {\cal T}'= \Big\{ \Omega=\cup_{\varepsilon\in I}\Omega_\varepsilon / \ \Omega_\varepsilon\in
{\cal T}'_\varepsilon\quad \forall \varepsilon\in I\Big\}$$
where ${\cal T}_\varepsilon$ and ${\cal T}'_\varepsilon$ are respectively the topology on $E_\varepsilon$ and $F_\varepsilon$ for all $\varepsilon\in I$.
We say that the diagonal topologies $ {\cal T}$ and ${\cal T}'$ are equivalent if for all $\varepsilon\in I$ the topologies ${\cal T}_\varepsilon$ and ${\cal T}'_\varepsilon$ are equivalent.
\end{defn}
\pesp

\begin{defn}\label{Def5}
A fractal family of topological spaces is a family

\begin{equation}\label{FF}
\Big(X_n^{j_n}, {\cal T}_n^{j_n}\Big)_{{{j_n}\atop n}{\in\atop\geq}{{\Lambda_n}\atop 0}}
\end{equation}
where

i) for all $n\geq0$, $\Lambda_n$ is an index set such that $\hbox{Card}\ (\Lambda_{n+1})>\hbox{Card}\ (\Lambda_n)$.

ii) for all $n\geq0$ and for all ${j_n}\in \Lambda_n$,  $\ (X_n^{j_n},\cal{T}_n^{j_n})$ is a topological space.

iii) for each $n\geq0$, the topologies $\cal{T}_n^{j_n}$ are equivalent for all ${j_n}\in \Lambda_n$.

iv) for all $n\geq0$ and for all $j_{n+1}\in \Lambda_{n+1}$, there exists a unique $j_n\in\Lambda_{n}$ such that
\begin{equation}
X_n^{j_n}\subset X_{n+1}^{j_{n+1}}\qquad  \hbox{and}\qquad {\cal T}_n^{j_n}= \Big\{O\cap X_n^{j_n}\ /\ O\in {\cal T}_{n+1}^{j_{n+1}}\Big\}.
\end{equation}

v) for all $n\geq0$, for all ${j_n}\in \Lambda_n$, there exists ${j_{n+1}}\in \Lambda_{n+1}$ such that
\begin{equation}\label{Topo}
{\cal T}_n^{j_n}\subset{\cal T}_{n+1}^{j_{n+1}}\qquad\hbox{and}\qquad{\cal T}_n^{j_n}= \Big\{O\cap X_n^{j_n}\ /\ O\in {\cal T}_{n+1}^{j_{n+1}}\Big\}.
\end{equation}
\end{defn}
\pesp

\begin{Exp}
Let us consider for example the index set $\Lambda_n=[1, 2^{n+1}]\cap\nN$ for all $n\geq0$. The number of topological spaces for successive iterations
will increase together with the index $n$, and we have:
\begin{itemize}
  \item for $n=0$, the family of topological spaces is given by  $(X_0^1,{\cal T}_0^1)$, $(X_0^2,{\cal T}_0^2)$.
  \item for $n=1$, the family of topological spaces is given by $(X_1^1,{\cal T}_1^1)$, $(X_1^2,{\cal T}_1^2)$, $(X_1^3,{\cal T}_1^3)$, $(X_1^4,{\cal T}_1^4)$.
  \item for $n=2$, the family of topological spaces is given by $(X_2^1,{\cal T}_2^1)$, $(X_2^2,{\cal T}_2^2)$, $(X_2^3,{\cal T}_2^3)$,
$(X_2^4,{\cal T}_2^4)$, $(X_2^5,{\cal T}_2^5)$, $(X_2^6,{\cal T}_2^6)$, $(X_2^7,{\cal T}_2^7)$, $(X_2^8,{\cal T}_2^8)$.
  \item more generally for $n$, the family is given by $2^{n+1}$ topological spaces.
\end{itemize}
\end{Exp}
\pesp

\begin{defn}
Let $\Big(X_n^{j_n}, {\cal T}_n^{j_n}\Big)_{{{j_n}\atop n}{\in\atop\geq}{{\Lambda_n}\atop 0}}$ be a fractal family of topological spaces. We call
the family $\Big({\cal T}_n^{j_n}\Big)_{{{j_n}\atop n}{\in\atop\geq}{{\Lambda_n}\atop 0}}$ a fractal topology.
\end{defn}
\pesp

\begin{prop}\label{prop2}
If $\di\Big(X_n^{j_n}, {\cal T}_n^{j_n}\Big)_{{{j_n}\atop n}{\in\atop\geq}{{\Lambda_n}\atop 0}}$ is a fractal family of topological spaces, then
for all $n>0$ and for all $j_0\in \Lambda_0$, there exist $j_1\in \Lambda_1, \ldots, j_n\in \Lambda_n$ such that
\begin{equation}
{\cal T}_0^{j_0}\subset {\cal T}_1^{j_1}\subset \ldots \subset {\cal T}_{n}^{j_n}.
\end{equation}
\end{prop}

\ni{\bf Proof.} For $n=1$, by property $v)$ of definition \ref{Def5},  there exists $j_1\in \Lambda_1$ such that using (\ref{Topo}) we have ${\cal T}_0^{j_0}\subset {\cal T}_1^{j_1}$.

By induction over $n>0$, suppose that there exist  $j_1\in \Lambda_1, \ldots, j_{n-1}\in \Lambda_{n-1}$ such that
${\cal T}_0^{j_0}\subset{\cal T}_1^{j_1}\subset   \ldots \subset {\cal T}_{n-1}^{j_{n-1}}$.
By property $v)$ of definition \ref{Def5},  there exists $j_n\in \Lambda_n$ such that we have ${\cal T}_{n-1}^{j_{n-1}}\subset {\cal T}_{n}^{j_n}$,
which completes the proof.

\rightline\Box

\begin{cor}
Under the condition of the proposition \ref{prop2}, the topology ${\cal T}_0^{j_0}$ is the weakest topology.
\end{cor}
\pesp

\begin{prop}
If $\di\Big(X_n^{j_n}, {\cal T}_n^{j_n}\Big)_{{{j_n}\atop n}{\in\atop\geq}{{\Lambda_n}\atop 0}}$ is a fractal family of topological spaces, then
for all $n>0$ and for all $j_n\in \Lambda_n$, there exist unique $j_0 \in  \Lambda_{0},\ldots, j_{n-1} \in  \Lambda_{n-1}$
such that
\begin{equation}\label{Sec}
X_0^{j_0}\subset  \ldots \subset X_{n-1}^{j_{n-1}}\subset X_{n}^{j_n}.
\end{equation}
\end{prop}

\ni{\bf Proof.}
For $i=1$ and $j_1\in \Lambda_{1}$, by definition \ref{Def5} iv), there exists a unique $j_0 \in \Lambda_{0}$ such that $ X_0^{j_0}\subset X_{1}^{j_{1}}$.

Let us suppose that for all $j_n\in \Lambda_{n}$ there exist unique $j_{n-1}\in \Lambda_{n-1},\ldots, j_0\in \Lambda_{0}$ such that
$X_{0}^{j_0}   \subset \ldots \subset X_{n-1}^{j_{n-1}}\subset  X_{n}^{j_{n}}$.

Let $j_{n+1}\ $ be in  $\ \Lambda_{n+1}$. By definition \ref{Def5} iv), there exists a unique $j_n\in \Lambda_{n}$ such that  $X_{n}^{j_{n}}\subset  X_{n+1}^{j_{n+1}}$.  By induction
there exist unique $j_{n-1}\in \Lambda_{n-1},\ldots, j_0\in \Lambda_{0}$ such that
$$X_{0}^{j_0}   \subset \ldots \subset X_{n-1}^{j_{n-1}}\subset  X_{n}^{j_{n}},$$ then we have
$X_{0}^{j_0}  \subset \ldots \subset X_{n-1}^{j_{n-1}}  \subset X_{n}^{j_{n}}\subset X_{n+1}^{j_{n+1}}$,
which gives the result.

\rightline\Box

\begin{prop}\label{prop4}
If $\di\Big(X_n^{j_n}, {\cal T}_n^{j_n}\Big)_{{{j_n}\atop n}{\in\atop\geq}{{\Lambda_n}\atop 0}}$ is a fractal family of topological spaces, then

i) for all $n\geq0$, $j_n\in \Lambda_{n}$ and $i>n$,  there exist $j_{n+1}\in \Lambda_{n+1},\ldots,j_i\in \Lambda_{i} $ such that the topology
${\cal T}_{n}^{j_n}$ is given by
\begin{equation}\label{F0i}
{\cal T}_{n}^{j_n}=\Big\{ O\cap X_{i-1}^{j_{i-1}}\cap \ldots\cap X_{n}^{j_n} \ / \  O\in{\cal T}_{i}^{j_{i}}  \Big\}.
\end{equation}

ii) for all $n\geq0$, $j_n\in \Lambda_{n}$ and $i<n$, there exist unique $j_{n-1}\in \Lambda_{n-1},\ldots,j_i\in \Lambda_{i} $ such that
\begin{equation}\label{F0ii}
{\cal T}_{i}^{j_i}=\Big\{ O\cap X_{n-1}^{j_{n-1}}\cap \ldots\cap X_{i}^{j_i} \ / \  O\in{\cal T}_{n}^{j_{n}}  \Big\}.
\end{equation}
\end{prop}

\ni{\bf Proof.}

i) Let $n\geq0$ and $j_n\in\Lambda_{n}$. Using definition \ref{Def5} v), there exists $j_{n+1}\in\Lambda_{n+1}$ such that
\begin{equation}\label{F1}
{\cal T}_{n}^{j_n}=\Big\{ O_{n+1}\cap X_{n}^{j_n} \ / \  O_{n+1}\in{\cal T}_{n+1}^{j_{n+1}}  \Big\}.
\end{equation}
By the same, for $j_{n+1}\in\Lambda_{n+1}$, there exists $j_{n+2}\in\Lambda_{n+2}$ such that
\begin{equation}\label{F2}
{\cal T}_{n+1}^{j_{n+1}}=\Big\{ O_{n+2}\cap X_{n+1}^{j_{n+1}} \ / \  O_{n+2}\in{\cal T}_{n+2}^{j_{n+2}}  \Big\}.
\end{equation}

Since in (\ref{F1})\  $O_{n+1}\in{\cal T}_{n+1}^{j_{n+1}}$, then using (\ref{F2}), there exists $\ O_{n+2}\in{\cal T}_{n+2}^{j_{n+2}}\ $ such that
$O_{n+1}=O_{n+2}\cap X_{n+1}^{j_{n+1}} $.
Therefore (\ref{F1}) becomes
\begin{equation}\label{F3}
{\cal T}_{n}^{j_n}=\Big\{ O_{n+2}\cap X_{n+1}^{j_{n+1}}\cap X_{n}^{j_n} \ / \  O_{n+2}\in{\cal T}_{n+2}^{j_{n+2}}  \Big\}.
\end{equation}

By induction over $i>n$, suppose that (\ref{F0i}) is true for $i=N$, with $N>n$, that is to say
 there exist $j_{n+1}\in \Lambda_{n+1},\ldots,j_{N}\in \Lambda_{N} $  such that
\begin{equation}\label{F5}
{\cal T}_{n}^{j_n}=\Big\{ O_N\cap X_{N-1}^{j_{N-1}}\cap \ldots\cap X_{n}^{j_n} \ / \  O_N\in{\cal T}_{N}^{j_{N}}  \Big\}.
\end{equation}

Let us prove (\ref{F0i}) for $i=N+1$. For $j_{N}\in \Lambda_{N}$, using definition \ref{Def5} v), there exists $j_{N+1}\in \Lambda_{N+1}$ such that
\begin{equation}\label{F6}
{\cal T}_{N}^{j_N}=\Big\{ O_{N+1}\cap X_{N}^{j_N} \ / \  O_{N+1}\in{\cal T}_{N+1}^{j_{N+1}}  \Big\}.
\end{equation}

Since in (\ref{F5})  $\ O_{N}\in{\cal T}_{N}^{j_{N}}$, then using (\ref{F6}), there exists $\ O_{N+1}\in{\cal T}_{N+1}^{j_{N+1}}\ $ such that
$O_{N}=O_{N+1}\cap X_{N}^{j_{N}} $. Then we obtain

$${\cal T}_{n}^{j_n}=\Big\{ O_{N+1}\cap X_{N}^{j_{N}}\cap X_{N-1}^{j_{N-1}}\cap \ldots\cap X_{n}^{j_n} \ / \  O_{N+1}\in{\cal T}_{N+1}^{j_{N+1}}  \Big\},$$
then (\ref{F0i}) is true for $i=N+1$.
\pesp

ii) Let $n\geq0$ and $j_n\in\Lambda_{n}$. By definition \ref{Def5} vi), there exists  a unique $j_{n-1}\in\Lambda_{n-1}$ such that
\begin{equation}\label{G1}
{\cal T}_{n-1}^{j_{n-1}}=\Big\{ O_n\cap X_{n-1}^{j_{n-1}} \ / \  O_n\in{\cal T}_{n}^{j_{n}}  \Big\}.
\end{equation}
By the same there exists a unique $j_{n-2}\in\Lambda_{n-2}$ such that
\begin{equation}\label{G2}
{\cal T}_{n-2}^{j_{n-2}}=\Big\{ O_{n-1}\cap X_{n-2}^{j_{n-2}} \ / \  O_{n-1}\in{\cal T}_{n-1}^{j_{n-1}}  \Big\}.
\end{equation}
Since in (\ref{G2}) $\ O_{n-1}\in{\cal T}_{n-1}^{j_{n-1}}\ $, then using (\ref{G1}) there exists $\ O_n\in{\cal T}_{n}^{j_{n}}\ $ such that
$$O_{n-1}=O_n\cap X_{n-1}^{j_{n-1}}.$$

Then (\ref{G2}) becomes
\begin{equation}\label{G3}
{\cal T}_{n-2}^{j_{n-2}}=\Big\{ O_{n}\cap X_{n-1}^{j_{n-1}}\cap X_{n-2}^{j_{n-2}} \ / \  O_{n}\in{\cal T}_{n}^{j_{n}}  \Big\}.
\end{equation}

By decreasing induction over $i<n$, suppose that (\ref{F0ii}) is true for $i=N$, with $N<n$, that is to say
there exist unique  $j_{n-1}\in \Lambda_{n-1},\ldots,j_N\in \Lambda_{N}$ such that
\begin{equation}\label{G4}
{\cal T}_{N}^{j_N}=\Big\{ O_n\cap X_{n-1}^{j_{n-1}}\cap \ldots\cap X_{N}^{j_N} \ / \  O_n\in{\cal T}_{n}^{j_{n}}  \Big\}.
\end{equation}
Let us prove (\ref{F0ii}) for $i=N-1$. Using definition \ref{Def5} iv), for $j_N\in \Lambda_{N}$, there exists a unique $j_{N-1}\in \Lambda_{N-1}$ such that
\begin{equation}\label{G5}
{\cal T}_{N-1}^{j_{N-1}}=\Big\{ O_N\cap X_{N-1}^{j_{N-1}}\ / \  O_N\in{\cal T}_{N}^{j_{N}}  \Big\}.
\end{equation}
Since in (\ref{G5}) $\ O_N\in{\cal T}_{N}^{j_{N}}\ $, then using (\ref{G4}), there exists $\ O_n\in{\cal T}_{n}^{j_{n}}\ $ such that
$$O_N=O_n\cap X_{n-1}^{j_{n-1}}\cap \ldots\cap X_{N}^{j_N}.$$
Therefore (\ref{G5}) becomes
$${\cal T}_{N-1}^{j_{N-1}}=\Big\{ O_n\cap X_{n-1}^{j_{n-1}}\cap \ldots\cap X_{N}^{j_N} \cap X_{N-1}^{j_{N-1}}\ / \  O_n\in{\cal T}_{n}^{j_{n}}  \Big\}$$
which gives (\ref{F0ii}) for $i=N-1$.

\rightline\Box

\begin{rem}\label{R2}

a) Under the conditions of the proposition \ref{prop4} i), we have $${\cal T}_{n}^{j_n}\subset{\cal T}_{n+1}^{j_{n+1}}\subset\ldots\subset{\cal T}_{i}^{j_i}.$$

b) Under the conditions of the proposition \ref{prop4} ii), we have $${\cal T}_{i}^{j_i}\subset\ldots\subset{\cal T}_{n-1}^{j_{n-1}}\subset{\cal T}_{n}^{j_n}.$$
\end{rem}
\pesp

\begin{cor}
Let $\di\Big(X_n^{j_n}, {\cal T}_n^{j_n}\Big)_{{{j_n}\atop n}{\in\atop\geq}{{\Lambda_n}\atop 0}}$ be a fractal family of topological spaces.

i) For all $n\geq0$, for all $j_n\in \Lambda_n$ and for all $i>n$, there exists $j_i\in \Lambda_i$ such that
\begin{equation}
{\cal T}_{n}^{j_n}=\Big\{ O\cap X_{n}^{j_n} \ / \  O\in{\cal T}_{i}^{j_{i}}  \Big\}.
\end{equation}

ii) For all $n>0$, for all $j_n\in \Lambda_n$ and for all $i<n$, there exists a unique
$j_i \in \Lambda_{i}$ such that
\begin{equation}
{\cal T}_{i}^{j_i}=\Big\{ O\cap X_{i}^{j_i} \ / \  O\in{\cal T}_{n}^{j_{n}}  \Big\}.
\end{equation}
\end{cor}

\ni{\bf Proof.}
 i) Let $n\geq0$ and $j_n\in \Lambda_n$.  By applying proposition \ref{prop4} i), for all $i>n$ there exist $j_{n+1}\in \Lambda_{n+1},\ldots,j_i\in \Lambda_{i} $ such that
\begin{equation}\label{F8}
{\cal T}_{n}^{j_n}=\Big\{ O\cap X_{i-1}^{j_{i-1}}\cap \ldots\cap X_{n}^{j_n} \ / \  O\in{\cal T}_{i}^{j_{i}}  \Big\}.
\end{equation}
By remark \ref{R2} a), ${\cal T}_{n}^{j_n}\subset{\cal T}_{n+1}^{j_{n+1}}\subset\ldots\subset{\cal T}_{i}^{j_i}$, which induces
$X_{n}^{j_n}\subset X_{n+1}^{j_{n+1}}\subset\ldots\subset X_{i}^{j_i}$, and then
$$X_{i-1}^{j_{i-1}}\cap \ldots\cap X_{n}^{j_n}=X_{n}^{j_n}.$$ Thus (\ref{F8}) becomes
\quad ${\cal T}_{n}^{j_n}=\Big\{ O\cap X_{n}^{j_n} \ / \  O\in{\cal T}_{i}^{j_{i}}  \Big\}.$
\pesp

ii) Let $n\geq0$ and $j_n\in \Lambda_n$. By applying proposition \ref{prop4} ii), for all $i<n$, there exist unique $j_{n-1}\in \Lambda_{n-1},\ldots,j_i\in \Lambda_{i} $
such that
\begin{equation}\label{F9}
{\cal T}_{i}^{j_i}=\Big\{ O\cap X_{n-1}^{j_{n-1}}\cap \ldots\cap X_{i}^{j_i} \ / \  O\in{\cal T}_{n}^{j_{n}}  \Big\}.
\end{equation}
By remark \ref{R2} b), ${\cal T}_{i}^{j_i}\subset\ldots\subset{\cal T}_{n-1}^{j_{n-1}}\subset{\cal T}_{n}^{j_n}$, which induces
$X_{i}^{j_i}\subset\ldots\subset X^{j_{n-1}}\subset X^{j_n}$, and then
$$X_{n-1}^{j_{n-1}}\cap \ldots\cap X_{i}^{j_i} =X_{i}^{j_i}.$$
Thus  (\ref{F9}) becomes \quad ${\cal T}_{i}^{j_i}=\Big\{ O\cap X_{i}^{j_i} \ / \  O\in{\cal T}_{n}^{j_{n}}  \Big\}.$

\rightline\Box

\begin{cor}\label{Co2}
Let $\di\Big(X_n^{j_n}, {\cal T}_n^{j_n}\Big)_{{{j_n}\atop n}{\in\atop\geq}{{\Lambda_n}\atop 0}}$ be a fractal family of topological spaces.

i) For all $j_0\in \Lambda_0$, for all $n>0$, there exists $j_n\in \Lambda_n$ such that
\begin{equation}
{\cal T}_{0}^{j_0}=\Big\{ O\cap X_{0}^{j_0} \ / \  O\in{\cal T}_{n}^{j_{n}}  \Big\}.
\end{equation}

ii) For all $n>0$, for all $j_n\in \Lambda_n$, there exists a unique
$j_0 \in \Lambda_{0}$ such that
\begin{equation}
{\cal T}_{0}^{j_0}=\Big\{ O\cap X_{0}^{j_0} \ / \  O\in{\cal T}_{n}^{j_{n}}  \Big\}.
\end{equation}
\end{cor}

\ni{\bf Proof.} i) We obtain the result by applying corollary \ref{Co2} i) for $n=0$ and $i=n$.

ii) We obtain the result by corollary \ref{Co2} ii) for $i=0$.

\rightline\Box

\begin{rem}
The family $\Big({\cal T}_n^{j_n}\Big)_{{{j_n}\atop n}{\in\atop\geq}{{\Lambda_n}\atop 0}}$ of topology
${\cal T}_n^{j_n}=\Big\{O\cap X_n^{j_n}\ /\ O\in {\cal T}_{n+1}^{j_{n+1}}\Big\}$  is called fractal topology because of the following:\pesp

1) there is a self-similarity in the construction of the induced topology for each $n\geq0$. The topology for each value of $n$ is obtained via the same process over bigger $n$.
The subspace topology that $X_0^{j_0}$ inherits from $X_1^{j_1}$ is the same than the one it inherits from $X_2^{j_2}$, etc.

2) the more $n$ increases, the more the topology ${\cal T}_n^{j_n}$ is strong, because of the appearance of new structures for each $n$.
\end{rem}
\pesp

\begin{defn}\label{Def4}
We call fractal topological space a family of sets $\di\Big(X_n^{j_n}\Big)_{{{j_n}\atop n}{\in\atop\geq}{{\Lambda_n}\atop 0}}$ endowed with a  fractal topology
$\di\Big({\cal T}_n^{j_n}\Big)_{{{j_n}\atop n}{\in\atop\geq}{{\Lambda_n}\atop 0}}$ such that
$\di\Big(X_n^{j_n}, {\cal T}_n^{j_n}\Big)_{{{j_n}\atop n}{\in\atop\geq}{{\Lambda_n}\atop 0}}$ is a fractal family of topological spaces.
\end{defn}
\pesp

\subsection{Fractal Manifold: a Fractal Topological Space}

Let us now examine the topological nature of a fractal manifold.
Let ${\cal M}$ be a fractal manifold. To prove that ${\cal M}$ is locally homeomorphic to a fractal topological space, we need to determine the family of
topologies associated to the family of sets
$\bigcup _{{\delta_0}...{\delta_{n}}}N_{\delta_0...\delta_{n}}^{\sigma_0...\sigma_{n}}$ given by theorem \ref{Th1}.
\pesp

\begin{prop}\label{pr2}
For all $n\geq0$, $\delta_i\in {{\cal R}_i}$ and $\sigma_i=\pm$, for $i=1,...,n$,
the set $N_{\delta_0...\delta_{n}}^{\sigma_0...\sigma_{n}}$ is a  Hausdorff topological space, and
if ${\cal T}_{\delta_0...\delta_{n}}^{\sigma_0...\sigma_{n}}$ is its associated topology,
then the set $\bigcup _{{\delta_0}...{\delta_{n}}} N_{\delta_0...\delta_{n}}^{\sigma_0...\sigma_{n}} $  is a diagonal topological space for the diagonal topology
\begin{equation}\label{F7}
{\cal T}_n^{\sigma_0...\sigma_{n}}= \Big\{ \Omega=\cup_{{\delta_0}...\delta_{n}}\Omega_{\delta_0...\delta_{n}}^{\sigma_0...\sigma_{n}}
\big/ \ \  \Omega_{\delta_0...\delta_{n}}^{\sigma_0...\sigma_{n}}\ \in  {\cal T}_{\delta_0...\delta_{n}}^{\sigma_0...\sigma_{n}}\ \
\forall {\delta_0}\in {{\cal R}_0},... ,\forall {\delta_{n}}\in {{\cal R}_n} \Big\}.
\end{equation}
\end{prop}

\ni{\bf Proof.} For $n=0$ and $\sigma_0=\pm$, the set $N_{\delta_0}^{\sigma_0}=\prod_{i=1}^{3}\Gamma_{i{\delta_0}}^{\sigma_0}\times\{{\delta_0}\}$
is the product of three graphs of mean functions given by (\ref{mom}), then $N_{\delta_0}^{\sigma_0}$ is an Hausdorff topological space and therefore
$\bigcup _{{\delta_0}
\in {\cal R}_0} N_{\delta_0}^{\sigma_0}$ is a disjoint union of Hausdorff topological spaces. If we denote by ${\cal T}_{\delta_0}^{\sigma_0}$
the topology on $N_{\delta_0}^{\sigma_0}$
for all $\delta_0\in {{\cal R}_0}$, then by the Definition \ref{DefTop}, we can associate to $\bigcup _{{\delta_0}
\in {\cal R}_0} N_{\delta_0}^{\sigma_0}$ the diagonal topology ${\cal T}_0^{\sigma_0}$ given by:
\begin{equation}
{\cal T}_0^{\sigma_0}= \Big\{ \Omega=\cup_{{\delta_0}\in {\cal R}_0}\Omega_{\delta_0}^{\sigma_0} \quad
\big/ \  \Omega_{\delta_0}^{\sigma_0}\  \in {\cal T}_{\delta_0}^{\sigma_0}\ \  \forall {\delta_0} \in {\cal R}_0 \Big\},
\end{equation}

\ni which makes $\bigcup _{{\delta_0}
\in {\cal R}_0}  N_{\delta_0}^{\sigma_0}$ a diagonal topological space.\pesp

By induction over $n\geq0$, suppose that for $\sigma_0=\pm$, ..., $\sigma_{n-1}=\pm$, the space
$N_{\delta_0...\delta_{n-1}}^{\sigma_0...\sigma_{n-1}}$ is a Hausdorff topological space, and that
the space $\bigcup _{{\delta_0}...{\delta_{n-1}}}
N_{\delta_0...\delta_{n-1}}^{\sigma_0...\sigma_{n-1}}$  is a diagonal topological space for the diagonal topology given by Definition \ref{DefTop}:

$${\cal T}_{n-1}^{\sigma_0...\sigma_{n-1}}= \Big\{ \Omega=\cup_{{\delta_0}...\delta_{n-1}}\Omega_{\delta_0...\delta_{n-1}}^{\sigma_0...\sigma_{n-1}}\
\big/ \  \Omega_{\delta_0...\delta_{n-1}}^{\sigma_0...\sigma_{n-1}}\ \in  {\cal T}_{\delta_0...\delta_{n-1}}^{\sigma_0...\sigma_{n-1}}\ \
\forall {\delta_0}\in {{\cal R}_0},... ,\forall {\delta_{n-1}}\in {{\cal R}_{n-1}} \Big\},$$

\ni where ${\cal T}_{\delta_0...\delta_{n-1}}^{\sigma_0...\sigma_{n-1}}$ is the topology on $N_{\delta_0...\delta_{n-1}}^{\sigma_0...\sigma_{n-1}}$
for all $\delta_0\in {{\cal R}_0}$,...,$\delta_{n-1}\in {{\cal R}_{n-1}}$.

Using Theorem \ref{Th1}, there exist a family of local homeomorphisms $\varphi_k$ and a family of translations $T_k$ for
$k\in [2^{n},2^{n+1}-1]\cap \nN$ at the $step(n)$ such that we have $2^n$ diagrams given by:
\pesp

\unitlength=1.2cm
\begin{picture}(11,3.7)

%%%%%++++step n

\put(6.4,2.7){ $\bigcup_{{\delta_0}...{\delta_{n}}}N_{\delta_0...\delta_{n}}^{\sigma_0...\sigma_{n-1}+}$}

%%%%%-----step n

\put(6.4,0.7){ $\bigcup_{{\delta_0}...{\delta_{n}}}N_{\delta_0...\delta_{n}}^{\sigma_0...\sigma_{n-1}-}$}

% les vecteurs
\put(4.8,2.7){$\varphi_k$}
\put(4.5,2.1){\vector(3,1){1.7}}

\put(7.8,1.8){$T_k$}
\put(7.5,2.3){\vector(0,-2){1.3}}

\put(4.5,0.8){$T_k\circ \varphi_k$}
\put(4.5,1.5){\vector(3,-1){1.7}}

%step (n-1)

\put(1.3,1.7){ $\bigcup_{{\delta_0}...{\delta_{n-1}}}N_{\delta_0...\delta_{n-1}}^{\sigma_0...\sigma_{n-1}}$}
\thicklines
\end{picture}
%\vskip0.1cm

\ni  for $\sigma_0=\pm$, ...,$\sigma_{n-1=\pm}$.

\ni For $\sigma_{n}=\pm$, the set\quad
$\bigcup _{{\delta_0}...{\delta_{n}}}N_{\delta_0...\delta_{n}}^{\sigma_0...\sigma_{n}} =\bigcup _{{\delta_0}...{\delta_{n}}}
\prod_{i=1}^{3}\Gamma_{i{\delta_0}...\delta_{n}}^{\sigma_0...\sigma_{n}}\times\{{\delta_{n}}\}\times\{{\delta_{n-1}}\}\times...\times\{{\delta_0}\} $
obtained at the $step(n)$ is a disjoint union of Hausdorff topological spaces:
indeed the set $\prod_{i=1}^{3}\Gamma_{i{\delta_0}...\delta_{n}}^{\sigma_0...\sigma_{n}}$ is the product of graphs of the function (\ref{mom}), then
$N_{\delta_0...\delta_{n}}^{\sigma_0...\sigma_{n}}$ is a Hausdorff topological space.
Therefore $\bigcup _{{\delta_0}...{\delta_{n}}}N_{\delta_0...\delta_{n}}^{\sigma_0...\sigma_{n}}$ is a disjoint union of Hausdorff topological spaces. Using
Definition \ref{DefTop}, we can associate to $\bigcup _{{\delta_0}...{\delta_{n}}}N_{\delta_0...\delta_{n}}^{\sigma_0...\sigma_{n}}$
the diagonal topology  ${\cal T} _n^{\sigma_0...\sigma_{n}}$ given by

$${\cal T}_n^{\sigma_0...\sigma_{n}}= \Big\{ \Omega=\cup_{{\delta_0}...\delta_{n}}\Omega_{\delta_0...\delta_{n}}^{\sigma_0...\sigma_{n}}\  /  \
\Omega_{\delta_0...\delta_{n}}^{\sigma_0...\sigma_{n}}\ \in  {\cal T}_{\delta_0...\delta_{n}}^{\sigma_0...\sigma_{n}}\ \  \forall {\delta_0}\in
{{\cal R}_0},... ,\forall {\delta_{n}}\in
{{\cal R}_n} \Big\},$$
where ${\cal T}_{\delta_0...\delta_{n}}^{\sigma_0...\sigma_{n}}$ is the topology on $N_{\delta_0...\delta_{n}}^{\sigma_0...\sigma_{n}}$ for all $\delta_0\in {{\cal R}_0}$,...,
$\delta_{n}\in {{\cal R}_n}$. Which gives that $\bigcup _{{\delta_0}...{\delta_{n}}} N_{\delta_0...\delta_{n}}^{\sigma_0...\sigma_{n}} $
is a diagonal topological space for the diagonal topology (\ref{F7}).

\rightline\Box

\begin{prop}\label{pr5}
For a given $n\geq0$, the diagonal topologies ${\cal T}_n^{\sigma_0...\sigma_{n}}$ are equivalent for $\sigma_0=\pm$, ...,$\sigma_n=\pm$.
\end{prop}

\ni{\bf Proof.} We know by proposition \ref{pr2} that
$\Big (\bigcup _{{\delta_0}...{\delta_{n}}} N_{\delta_0...\delta_{n}}^{\sigma_0...\sigma_{n}},\ {\cal T}_n^{\sigma_0...\sigma_{n}}\Big)$
is a diagonal topological space for
$n\geq0$, $\sigma_0=\pm$, ..., $\sigma_{n}=\pm$. Let $n\geq0$ and $\delta_0\in {\cal R}_0,\ldots,\delta_n\in {\cal R}_n$.
Since  $\di N_{\delta_0...\delta_{n}}^{\sigma_0...\sigma_{n}}=
\prod_{i=1}^{3}\Gamma_{i{\delta_0}...\delta_{n}}^{\sigma_0...\sigma_{n}}\times\{{\delta_{n}}\}\times\{{\delta_{n-1}}\}\times...\times\{{\delta_0}\}$
is the product of three graphs of function given by (\ref{mom}), then for $\sigma_0=\pm,\ldots,\sigma_n=\pm$\quad
the spaces $N_{\delta_0...\delta_{n}}^{\sigma_0...\sigma_{n}}$ are homeomorphic to $\rR^3$, which means that for
$\sigma_0=\pm,\ldots,\sigma_n=\pm$\quad the topologies ${\cal T}_{\delta_0...\delta_{n}}^{\sigma_0...\sigma_{n}}$
are equivalent. By definition \ref{ET}, the diagonal topologies
${\cal T}_n^{\sigma_0...\sigma_{n}}$ on $\di\bigcup _{{\delta_0}...{\delta_{n}}} N_{\delta_0...\delta_{n}}^{\sigma_0...\sigma_{n}}$
are equivalent for $\sigma_0=\pm$, ...,$\sigma_n=\pm$.

\rightline\Box

\begin{prp}\label{pr1}
For all $n\geq0$, for all $\sigma_0=\pm,\ldots,\sigma_n=\pm$ and for all $\delta_0\in {\cal R}_0,\ldots,\delta_n\in{\cal R}_n$, the set
$N_{\delta_0...\delta_{n}}^{\sigma_0...\sigma_{n}}$ is identical to the set
$N_{\delta_0...\delta_{n}0}^{\sigma_0...\sigma_{n}\sigma_{n+1}}$ for $\sigma_{n+1}=\pm$.
\end{prp}

\ni{\bf Proof.}
Let $n\geq0$ and $\delta_0\in {\cal R}_0,\ldots,\delta_n\in {\cal R}_n$ and $\sigma_0=\pm,\ldots,\sigma_n=\pm$. By (\ref{F10}),

$$ N_{\delta_0...\delta_{n}}^{\sigma_0...\sigma_{n}}=\prod_{i=1}^3\Gamma_{i\delta_{0}...\delta_{n}}^{\sigma_0...\sigma_{n}}
\times\{\delta_{n}\}
\times...\times\{\delta_{0}\}$$
where for $i=1,2,3$, $\Gamma_{i\delta_{0}...\delta_{n}}^{\sigma_0...\sigma_{n}}$ is the graph of the function $F^{\sigma_0...\sigma_{n}}_{i\delta_0...\delta_{n}}$.

Since $\delta_{n+1}\in {\cal R}_{n+1}=[0, \varepsilon_{n+1}]$, then we have for $\sigma_{n+1}=\pm$
$$F^{\sigma_0...\sigma_{n+1}}_{i\delta_0...\delta_{n}0}(x)=\lim_{\delta_{n+1}\rightarrow 0} F^{\sigma_0...\sigma_{n+1}}_{i\delta_0...\delta_{n}\delta_{n+1}}(x)
=F^{\sigma_0...\sigma_{n}}_{i\delta_0...\delta_{n}}(x).$$

Therefore we have for $\sigma_{n+1}=\pm$ and $i=1,2,3$
$$\Gamma_{i\delta_{0}...\delta_{n}0}^{\sigma_0...\sigma_{n}\sigma_{n+1}}=\Gamma_{i\delta_{0}...\delta_{n}}^{\sigma_0...\sigma_{n}},$$
which yields that for $\sigma_{n+1}=\pm$
$$\prod_{i=1}^3\Gamma_{i\delta_{0}...\delta_{n}0}^{\sigma_0...\sigma_{n}\sigma_{n+1}}\times\{0\}
\times\{\delta_{n}\}
\times...\times\{\delta_{0}\}=\prod_{i=1}^3\Gamma_{i\delta_{0}...\delta_{n}}^{\sigma_0...\sigma_{n}}
\times\{\delta_{n}\}
\times...\times\{\delta_{0}\}$$
by identification of the points $(X,Y,Z,0,\delta_{n},\ldots,\delta_{0})$ and $(X,Y,Z,\delta_{n},\ldots,\delta_{0})$ where

\ni$X\in\Gamma_{1\delta_{0}...\delta_{n}}^{\sigma_0...\sigma_{n}}\ $,
$Y\in \Gamma_{2\delta_{0}...\delta_{n}}^{\sigma_0...\sigma_{n}}\ $ and
$\ Z\in \Gamma_{3\delta_{0}...\delta_{n}}^{\sigma_0...\sigma_{n}}$.
Thus we obtain

$$N_{\delta_0...\delta_{n}0}^{\sigma_0...\sigma_{n}\sigma_{n+1}}=N_{\delta_0...\delta_{n}}^{\sigma_0...\sigma_{n}}\quad \hbox{for} \quad \sigma_{n+1}=\pm.$$

\rightline\Box

\begin{lem} \label{L1}
For all $n\geq0$ and for all $\sigma_0=\pm,\ldots,\sigma_n=\pm$, the set $\di\bigcup _{{\delta_0}...{\delta_{n}}} N_{\delta_0...\delta_{n}}^{\sigma_0...\sigma_{n}}$
is identical to the subset \quad$\di\bigcup _{{\delta_0}...{\delta_{n}}} N_{\delta_0...\delta_{n}0}^{\sigma_0...\sigma_{n}\sigma_{n+1}}$
of \quad $\di\bigcup _{{\delta_0}...{\delta_{n+1}}} N_{\delta_0...\delta_{n+1}}^{\sigma_0...\sigma_{n+1}}\ $ for $\ \sigma_{n+1}=\pm$ .
\end{lem}

\ni{\bf Proof.}
Let $n\geq0$, $\sigma_0=\pm,\ldots,\sigma_{n}=\pm$. By property \ref{pr1}, for all $\delta_0\in {\cal R}_0,\ldots,\delta_n\in{\cal R}_n$, we have
$$N_{\delta_0...\delta_{n}}^{\sigma_0...\sigma_{n}}=N_{\delta_0...\delta_{n}0}^{\sigma_0...\sigma_{n}\sigma_{n+1}} \quad \hbox{for} \quad \sigma_{n+1}=\pm.$$
Then
\begin{equation}\label{F11}
\bigcup _{{\delta_0}...{\delta_{n}}} N_{\delta_0...\delta_{n}}^{\sigma_0...\sigma_{n}}=
\bigcup _{{\delta_0}...{\delta_{n}}} N_{\delta_0...\delta_{n}0}^{\sigma_0...\sigma_{n}\sigma_{n+1}} \quad \hbox{for} \quad \sigma_{n+1}=\pm.
\end{equation}
Since

$$\bigcup _{{\delta_0}...\delta_{n}{\delta_{n+1}}} N_{\delta_0...\delta_{n+1}}^{\sigma_0...\sigma_{n+1}}=
\Big(\bigcup _{{{\delta_0}\atop {\delta_{n+1}}}{...\atop\neq}{{\delta_{n}}\atop 0}}N_{\delta_0...\delta_{n+1}}^{\sigma_0...\sigma_{n+1}}\Big)
\bigcup \Big(\bigcup _{{{\delta_0}\atop {\delta_{n+1}}}{...\atop =}{{\delta_{n}}\atop 0}}N_{\delta_0...\delta_{n}0}^{\sigma_0...\sigma_{n}\sigma_{n+1}}\Big),$$

then
$\di\bigcup _{{\delta_0}...{\delta_{n}}} N_{\delta_0...\delta_{n}0}^{\sigma_0...\sigma_{n}\sigma_{n+1}} $
is the subset of $\di\bigcup _{{\delta_0}...{\delta_{n+1}}} N_{\delta_0...\delta_{n+1}}^{\sigma_0...\sigma_{n+1}}$ that corresponds to the value $\delta_{n+1}=0$.
We deduce from (\ref{F11}) that for $\sigma_{n+1}=\pm$
$$\bigcup _{{\delta_0}...{\delta_{n}}} N_{\delta_0...\delta_{n}}^{\sigma_0...\sigma_{n}}=
\bigcup _{{\delta_0}...{\delta_{n}}} N_{\delta_0...\delta_{n}0}^{\sigma_0...\sigma_{n}\sigma_{n+1}}\quad \hbox{in} \quad
\bigcup _{{\delta_0}...{\delta_{n+1}}} N_{\delta_0...\delta_{n+1}}^{\sigma_0...\sigma_{n+1}}.$$

\rightline\Box

\begin{rem}
The previous identification is true for $\sigma_{n+1}=\pm$, which means that the set
$\di\bigcup _{{\delta_0}...{\delta_{n}}} N_{\delta_0...\delta_{n}}^{\sigma_0...\sigma_{n}}$ is duplicated in two sets
$\di\bigcup _{{\delta_0}...{\delta_{n}}} N_{\delta_0...\delta_{n}0}^{\sigma_0...\sigma_{n}+}$ and
$\di\bigcup _{{\delta_0}...{\delta_{n}}} N_{\delta_0...\delta_{n}0}^{\sigma_0...\sigma_{n}-}$ respectively in
$\di\bigcup _{{\delta_0}...{\delta_{n+1}}} N_{\delta_0...\delta_{n+1}}^{\sigma_0...\sigma_{n}+}$ and
$\di\bigcup _{{\delta_0}...{\delta_{n+1}}} N_{\delta_0...\delta_{n+1}}^{\sigma_0...\sigma_{n}-}$.
\end{rem}
\pesp

\begin{thm}\label{Th6}
If for all $n\geq0$ and for all $\sigma_0=\pm,\ldots,\sigma_n=\pm$,
the set  $\di\bigcup _{{\delta_0}...{\delta_{n}}} N_{\delta_0...\delta_{n}}^{\sigma_0...\sigma_{n}}$ is endowed with the
diagonal topology ${\cal T}_{n}^{\sigma_0...\sigma_{n}}$ given by (\ref{F7}), then
\begin{equation}\label{IN}
{\cal T}_{n}^{\sigma_0...\sigma_{n}}=\Big\{O\cap \di\bigcup _{{\delta_0}...{\delta_{n}}}
N_{\delta_0...\delta_{n}}^{\sigma_0...\sigma_{n}}\ /\  O\in {\cal T}_{n+1}^{\sigma_0...\sigma_{n}\sigma_{n+1}}\Big\}\quad \hbox{for}\quad \sigma_{n+1}=\pm.
\end{equation}
\end{thm}

\ni{\bf Proof.} To prove (\ref{IN}), we have to prove that for $n\geq0$ and $\sigma_0=\pm,\ldots,\sigma_n=\pm$
\begin{equation}\label{IN1}
i) \qquad \qquad {\cal T}_{n}^{\sigma_0...\sigma_{n}}\subset\Big\{O\cap \di\bigcup _{{\delta_0}...{\delta_{n}}}
N_{\delta_0...\delta_{n}}^{\sigma_0...\sigma_{n}}\ /\  O\in {\cal T}_{n+1}^{\sigma_0...\sigma_{n}\sigma_{n+1}}\Big\}
\end{equation}

\begin{equation}\label{IN2}
ii) \qquad \qquad {\cal T}_{n}^{\sigma_0...\sigma_{n}}\supset\Big\{O\cap \di\bigcup _{{\delta_0}...{\delta_{n}}}
N_{\delta_0...\delta_{n}}^{\sigma_0...\sigma_{n}}\ /\  O\in {\cal T}_{n+1}^{\sigma_0...\sigma_{n}\sigma_{n+1}}\Big\}
\end{equation}
where ${\cal T}_{n}^{\sigma_0...\sigma_{n}}$ is the diagonal topology on
$\bigcup _{{\delta_0}...{\delta_{n}}} N_{\delta_0...\delta_{n}}^{\sigma_0...\sigma_{n}}$.
\pesp

i) Let $n\geq0$ and $\sigma_0=\pm,\ldots,\sigma_n=\pm$. Let $\Omega$ be an open set of
$\di\bigcup _{{\delta_0}...{\delta_{n}}} N_{\delta_0...\delta_{n}}^{\sigma_0...\sigma_{n}}$ for the diagonal topology
${\cal T}_{n}^{\sigma_0...\sigma_{n}}$. By the definition of
${\cal T}_{n}^{\sigma_0...\sigma_{n}}$, $\di\Omega=\bigcup _{{\delta_0}...{\delta_{n}}} \Omega_{\delta_0...\delta_{n}}^{\sigma_0...\sigma_{n}}$, where
 $\Omega_{\delta_0...\delta_{n}}^{\sigma_0...\sigma_{n}}$ is an open set of
$N_{\delta_0...\delta_{n}}^{\sigma_0...\sigma_{n}}$ for all
$\delta_0\in {\cal R}_0,\ldots,\delta_n\in{\cal R}_n$.

By property \ref{pr1}, for all $\delta_0\in {\cal R}_0,\ldots,\delta_n\in{\cal R}_n$ and for $\sigma_{n+1}=\pm$, we have
$$N_{\delta_0...\delta_{n}}^{\sigma_0...\sigma_{n}}=N_{\delta_0...\delta_{n}0}^{\sigma_0...\sigma_{n+1}}$$
then any open set of $N_{\delta_0...\delta_{n}}^{\sigma_0...\sigma_{n}}$ is an open set of $N_{\delta_0...\delta_{n}0}^{\sigma_0...\sigma_{n+1}}$ for $\sigma_{n+1}=\pm$.

Since for all $\ \delta_0\in {\cal R}_0,\ldots,\delta_n\in{\cal R}_n$, \ $\Omega_{\delta_0...\delta_{n}}^{\sigma_0...\sigma_{n}}$ is an open set of
$N_{\delta_0...\delta_{n}}^{\sigma_0...\sigma_{n}}$ then for all $\delta_0\in {\cal R}_0,\ldots,\delta_n\in{\cal R}_n$,
 $\Omega_{\delta_0...\delta_{n}}^{\sigma_0...\sigma_{n}}$ is also an open set of
$N_{\delta_0...\delta_{n}0}^{\sigma_0...\sigma_{n+1}}$ for $\sigma_{n+1}=\pm$.

Let us consider the set
\begin{equation}\label{O1}
O=\di\bigcup _{{\delta_0}...{\delta_n}{\delta_{n+1}}} O_{\delta_0...\delta_{n}\delta_{n+1}}^{\sigma_0...\sigma_{n}\sigma_{n+1}}
\end{equation}
defined by

\begin{equation}\label{O2}
O_{\delta_0...\delta_{n}\delta_{n+1}}^{\sigma_0...\sigma_{n}\sigma_{n+1}}=
\left\{
  \begin{array}{ll}
    \Omega_{\delta_0...\delta_{n}}^{\sigma_0...\sigma_{n}}& \forall \delta_0,\ldots,\delta_n\quad  \hbox{and for} \quad \delta_{n+1}=0\\
    \emptyset & \forall \delta_0,\ldots,\delta_n\quad  \hbox{and for} \quad \delta_{n+1}\neq0.\\
  \end{array}
\right.
\end{equation}

We first have to verify that $\ O\in {\cal T}_{n+1}^{\sigma_0...\sigma_{n}\sigma_{n+1}}$:
\quad for all $\delta_0,\ldots,\delta_n$, and for $\delta_{n+1}=0$,
$\ O_{\delta_0...\delta_{n}0}^{\sigma_0...\sigma_{n}\sigma_{n+1}}=\Omega_{\delta_0...\delta_{n}}^{\sigma_0...\sigma_{n}}\ $
is an open set of
$\ N_{\delta_0...\delta_{n}0}^{\sigma_0...\sigma_{n}\sigma_{n+1}}$. For all $\delta_0,\ldots,\delta_n$, and for $\delta_{n+1}\neq0$, $\emptyset$ is an open set of
$N_{\delta_0...\delta_{n}{\delta_{n+1}}}^{\sigma_0...\sigma_{n}\sigma_{n+1}}$, then  $O\in {\cal T}_{n+1}^{\sigma_0...\sigma_{n}\sigma_{n+1}}$.

Secondly we have to verify that $O\cap \di\bigcup _{{\delta_0}...{\delta_{n}}}
N_{\delta_0...\delta_{n}}^{\sigma_0...\sigma_{n}}=\Omega$:\quad
indeed by (\ref{O1}),

$$O\cap \bigcup _{{\delta_0}...{\delta_{n}}}N_{\delta_0...\delta_{n}}^{\sigma_0...\sigma_{n}}=
\Big(\bigcup _{{\delta_0}...{\delta_n}{\delta_{n+1}}} O_{\delta_0...\delta_{n}\delta_{n+1}}^{\sigma_0...\sigma_{n}\sigma_{n+1}}\Big)
\cap \Big(\bigcup _{{\delta_0}...{\delta_{n}}}N_{\delta_0...\delta_{n}}^{\sigma_0...\sigma_{n}}\Big)$$

$$=\Big((\bigcup _{\delta_{n+1}\neq0}\bigcup _{{\delta_0}...{\delta_{n}}}O_{\delta_0...\delta_{n}\delta_{n+1}}^{\sigma_0...\sigma_{n}\sigma_{n+1}})\cup
(\bigcup _{\delta_{n+1}=0}\bigcup _{{\delta_0}...{\delta_{n}}}O_{\delta_0...\delta_{n}\delta_{n+1}}^{\sigma_0...\sigma_{n}\sigma_{n+1}})\Big)
\cap \Big(\bigcup _{{\delta_0}...{\delta_{n}}}N_{\delta_0...\delta_{n}}^{\sigma_0...\sigma_{n}}\Big).$$

By (\ref{O2}), we obtain

$$O\cap \bigcup _{{\delta_0}...{\delta_{n}}}N_{\delta_0...\delta_{n}}^{\sigma_0...\sigma_{n}}= \Big(\emptyset\cup
(\bigcup _{{\delta_0}...{\delta_{n}}}\Omega_{\delta_0...\delta_{n}}^{\sigma_0...\sigma_{n}})\Big)
\cap \Big(\bigcup _{{\delta_0}...{\delta_{n}}}N_{\delta_0...\delta_{n}}^{\sigma_0...\sigma_{n}}\Big)$$
$$= \Big(\bigcup _{{\delta_0}...{\delta_{n}}}\Omega_{\delta_0...\delta_{n}}^{\sigma_0...\sigma_{n}}\Big)
\cap \Big(\bigcup _{{\delta_0}...{\delta_{n}}}N_{\delta_0...\delta_{n}}^{\sigma_0...\sigma_{n}}\Big)$$
$$=\Omega\cap \Big(\bigcup _{{\delta_0}...{\delta_{n}}}N_{\delta_0...\delta_{n}}^{\sigma_0...\sigma_{n}}\Big)=\Omega,$$
which induces the inclusion
$${\cal T}_{n}^{\sigma_0...\sigma_{n}}\subset\Big\{O\cap \di\bigcup _{{\delta_0}...{\delta_{n}}}
N_{\delta_0...\delta_{n}}^{\sigma_0...\sigma_{n}}\ /\  O\in {\cal T}_{n+1}^{\sigma_0...\sigma_{n}\sigma_{n+1}}\Big\}.$$
\pesp

ii) Inversely, let $n\geq0$, $\sigma_0=\pm,\ \ldots,\ \sigma_{n+1}=\pm$, and let us consider $\Omega$ in the set $\Big\{O\cap \di\bigcup _{{\delta_0}...{\delta_{n}}}
N_{\delta_0...\delta_{n}}^{\sigma_0...\sigma_{n}}\ /\  O\in {\cal T}_{n+1}^{\sigma_0...\sigma_{n}\sigma_{n+1}}\Big\}$, then there exists
$O\in {\cal T}_{n+1}^{\sigma_0...\sigma_{n}\sigma_{n+1}}$
such that
\begin{equation}\label{F4}
\Omega=O\cap \di\bigcup _{{\delta_0}...{\delta_{n}}}
N_{\delta_0...\delta_{n}}^{\sigma_0...\sigma_{n}}.
\end{equation}
Since $\ {\cal T}_{n+1}^{\sigma_0...\sigma_{n}\sigma_{n+1}}\ $ is the diagonal topology on
$\ \di\bigcup _{{\delta_0}...{\delta_{n+1}}} N_{\delta_0...\delta_{n+1}}^{\sigma_0...\sigma_{n+1}}\ $ given by (\ref{F7}), then for all
$\ \delta_0\in{\cal R}_0,\ldots,\delta_{n+1}\in{\cal R}_{n+1},\ $ there exists
$\ O_{\delta_0...\delta_{n+1}}^{\sigma_0...\sigma_{n+1}}\in {\cal T}_{{\delta_0}...\delta_{n+1}}^{\sigma_0...\sigma_{n+1}} $ such that
$O=\bigcup_{{\delta_0}...\delta_{n+1}}O_{\delta_0...\delta_{n+1}}^{\sigma_0...\sigma_{n+1}}$.\pesp

By substitution in (\ref{F4}), we obtain
$$\Omega=\Big(\bigcup_{{\delta_0}...\delta_{n+1}}O_{\delta_0...\delta_{n+1}}^{\sigma_0...\sigma_{n+1}}\Big)\cap \Big(\bigcup _{{\delta_0}...{\delta_{n}}}
N_{\delta_0...\delta_{n}}^{\sigma_0...\sigma_{n}}\Big)$$
$$=\Big((\bigcup _{\delta_{n+1}\neq0}\bigcup _{{\delta_0}...{\delta_{n}}}O_{\delta_0...\delta_{n}\delta_{n+1}}^{\sigma_0...\sigma_{n+1}})\cup
(\bigcup _{\delta_{n+1}=0}\bigcup _{{\delta_0}...{\delta_{n}}}O_{\delta_0...\delta_{n}\delta_{n+1}}^{\sigma_0...\sigma_{n+1}})\Big)
\cap \Big(\bigcup _{{\delta_0}...{\delta_{n}}} N_{\delta_0...\delta_{n}}^{\sigma_0...\sigma_{n}}\Big).$$
Since by lemma \ref{L1},
$$\bigcup _{{\delta_0}...{\delta_{n}}} N_{\delta_0...\delta_{n}}^{\sigma_0...\sigma_{n}}=
\bigcup _{{\delta_0}...{\delta_{n}}} N_{\delta_0...\delta_{n}0}^{\sigma_0...\sigma_{n+1}},$$
then

$$\Omega=\Big((\bigcup _{\delta_{n+1}\neq0}\bigcup _{{\delta_0}...{\delta_{n}}}O_{\delta_0...\delta_{n}\delta_{n+1}}^{\sigma_0...\sigma_{n+1}})\cup
(\bigcup _{\delta_{n+1}=0}\bigcup _{{\delta_0}...{\delta_{n}}}O_{\delta_0...\delta_{n}\delta_{n+1}}^{\sigma_0...\sigma_{n+1}})\Big)
\cap \Big(\bigcup _{{\delta_0}...{\delta_{n}}} N_{\delta_0...\delta_{n}0}^{\sigma_0...\sigma_{n+1}}\Big)$$

$$=\Big(( \bigcup _{\delta_{n+1}\neq0}\bigcup _{{\delta_0}...{\delta_{n}}}O_{\delta_0...\delta_{n}\delta_{n+1}}^{\sigma_0...\sigma_{n+1}})\cap
\bigcup _{{\delta_0}...{\delta_{n}}}N_{\delta_0...\delta_{n}0}^{\sigma_0...\sigma_{n+1}}\Big)
\cup\Big((\bigcup _{\delta_{n+1}=0}\bigcup _{{\delta_0}...{\delta_{n}}}O_{\delta_0...\delta_{n}\delta_{n+1}}^{\sigma_0...\sigma_{n+1}})\cap
\bigcup _{{\delta_0}...{\delta_{n}}}N_{\delta_0...\delta_{n}0}^{\sigma_0...\sigma_{n+1}}\Big).$$

Since
$$\Big( \di \bigcup _{\delta_{n+1}\neq0}\bigcup _{{\delta_0}...{\delta_{n}}}O_{\delta_0...\delta_{n}\delta_{n+1}}^{\sigma_0...\sigma_{n+1}}\Big)\cap
\Big(\bigcup _{{\delta_0}...{\delta_{n}}}N_{\delta_0...\delta_{n}0}^{\sigma_0...\sigma_{n+1}}\Big)=\emptyset,$$
then

$$\Omega=\Big(\bigcup _{\delta_{n+1}=0}\bigcup _{{\delta_0}...{\delta_{n}}}O_{\delta_0...\delta_{n}\delta_{n+1}}^{\sigma_0...\sigma_{n+1}}\Big)\cap
\Big(\bigcup _{{\delta_0}...{\delta_{n}}}N_{\delta_0...\delta_{n}0}^{\sigma_0...\sigma_{n+1}}\Big)$$

$$=\Big(\bigcup _{{\delta_0}...{\delta_{n}}}O_{\delta_0...\delta_{n}0}^{\sigma_0...\sigma_{n+1}}\Big)\cap
\Big(\bigcup _{{\delta_0}...{\delta_{n}}}N_{\delta_0...\delta_{n}0}^{\sigma_0...\sigma_{n+1}}\Big)$$

Since for all $\delta_0\in{\cal R}_0,\ldots,\delta_{n}\in{\cal R}_{n}$,  $ \di O_{{\delta_0}...{\delta_n}{0}}^{\sigma_0...\sigma_{n+1}}$ is an open set of
$N_{\delta_0...\delta_{n}0}^{\sigma_0...\sigma_{n}\sigma_{n+1}}$, then for all $\delta_0\in{\cal R}_0,\ldots,\delta_{n}\in{\cal R}_{n}$,
$ \di O_{{\delta_0}...{\delta_n}{0}}^{\sigma_0...\sigma_{n+1}}$ is an open set of $\di N_{\delta_0...\delta_{n}}^{\sigma_0...\sigma_{n}}$
by property \ref{pr1}, and then
$\di \bigcup _{{\delta_0}...{\delta_{n}}} \Omega_{{\delta_0}...{\delta_n}{0}}^{\sigma_0...\sigma_{n+1}}$ is an open set of
$\di \bigcup _{{\delta_0}...{\delta_{n}}}N_{\delta_0...\delta_{n}}^{\sigma_0...\sigma_{n}}$.

Therefore we obtain
$\di\Omega=\bigcup _{{\delta_0}...{\delta_{n}}} O_{{\delta_0}...{\delta_n}{0}}^{\sigma_0...\sigma_{n+1}}\in {\cal T}_{n}^{\sigma_0...\sigma_{n}}$
which proves (\ref{IN2}).

\rightline\Box

We introduce the following index set:\pesp

\ni{\bf Notation:}
We denote for all $n\geq0$,
\begin{equation}\label{DefSet}
\Lambda_n=\{\sigma_0\ldots\sigma_{n}\ /\ \sigma_0=\pm,\ldots,\sigma_n=\pm\}
\end{equation}
 where the cardinal of $\Lambda_n$ is $2^{n+1}$, that is to say:

$\Lambda_0=\{\sigma_0\ /\ \sigma_0=\pm\}=\{+, -\}$ with cardinal $2$,

$\Lambda_1=\{\sigma_0\sigma_{1}\ /\ \sigma_0=\pm,\sigma_1=\pm\}=\{++,+-,-+,--\}$ with cardinal $2^2$,

$\Lambda_2=\{\sigma_0\sigma_{1}\sigma_{2}\ /\ \sigma_0=\pm,\sigma_1=\pm,\sigma_2=\pm\}$

$\quad=\{+++,++-,+-+,+--,-++,-+-,--+,---\}$ with cardinal $2^3$.
\pesp

\begin{thm}\label{Th5}
The family of diagonal topological spaces
\begin{equation}\label{Ft1}
\Big(\bigcup _{{\delta_0}...{\delta_{n}}} N_{\delta_0...\delta_{n}}^{j_n},\ {\cal T}_{n}^{j_n}\Big)_
{{j_n\atop n}{\in\Lambda_n\atop\geq0}}
\end{equation}
where $\Lambda_n$ is given by (\ref{DefSet}) for $n\geq0$,
is a fractal family of topological spaces.
\end{thm}

\ni{\bf Proof.} We have to prove that the family (\ref{Ft1}) verifies the properties from i) to v) of the definition \ref{Def5}.
\pesp

i) The cardinal of the index set $\Lambda_n$ is $2^{n+1}$ and the cardinal of $\Lambda_{n+1}$ is $2^{n+2}$, then the cardinal of $\Lambda_{n+1}$
is strictly greater than the cardinal of $\Lambda_n$ for all $n\geq0$.
\pesp

ii) For all $n\geq0$ and for all $j_n\in\Lambda_n$, the space
$\Big( \bigcup _{{\delta_0}...{\delta_{n}}} N_{\delta_0...\delta_{n}}^{j_n},\ {\cal T}_{n}^{j_n}\Big)$ is a diagonal topological
space by proposition \ref{pr2}, then it is a topological space.
\pesp

iii) For each $n\geq0$, the topologies ${\cal T}_{n}^{j_n}$ are equivalent for all $j_n\in \Lambda_n$ by proposition \ref{pr5}.
\pesp

iv) Let $n\geq0$ and $j_{n+1}\in\Lambda_{n+1}$. By (\ref{DefSet}), there exist $\sigma_0=\pm,\ldots,\sigma_{n}=\pm,\sigma_{n+1}=\pm$ such that
 $\ j_{n+1}=\sigma_0\ldots\sigma_{n}\sigma_{n+1}$.
By lemma \ref{L1}, the unique $\ j_n\in \Lambda_n\ $ such that $$\bigcup _{{\delta_0}...{\delta_{n}}} N_{\delta_0...\delta_{n}}^{j_n}\subset
\bigcup _{{\delta_0}...{\delta_{n+1}}} N_{\delta_0...\delta_{n+1}}^{j_{n+1}}$$
is given
by $j_n=\sigma_0\ldots\sigma_{n}$, and ${\cal T}_{n}^{j_n}$ is the induced topology given by (\ref{IN}).
\pesp

v) To prove that the condition v) of the definition \ref{Def5} is satisfied by the family (\ref{Ft1}), we have to prove that
${\cal T}_{n}^{j_n}\subset{\cal T}_{n+1}^{j_{n+1}}$ for all $n\geq0$.

Let $n\geq0$ and $j_{n}\in\Lambda_{n}$. By (\ref{DefSet}), there exist $\sigma_0=\pm,\ldots,\sigma_{n}=\pm$ such that
 $j_{n}=\sigma_0\ldots\sigma_{n}$. By theorem \ref{Th6}, for $\sigma_{n+1}=\pm$ and $j_{n+1}=\sigma_0\ldots\sigma_{n}\sigma_{n+1}\in \Lambda_{n+1}$,
we have
\begin{equation}
{\cal T}_{n}^{j_n}=\Big\{O\cap\Big(\bigcup _{{\delta_0}...{\delta_{n}}} N_{\delta_0...\delta_{n}}^{j_n}\Big)\quad /\quad O\in {\cal T}_{n+1}^{j_{n+1}}\Big\}.
\end{equation}

Let us consider $\Omega\in {\cal T}_{n}^{j_n}$, then $\Omega=\bigcup _{{\delta_0}...{\delta_{n}}}\Omega_{\delta_0...\delta_{n}}^{j_n}$, where
$\Omega_{\delta_0...\delta_{n}}^{j_n}$ is an open set of $N_{\delta_0...\delta_{n}}^{j_n}$ for all $\delta_0...\delta_{n}$.

By theorem \ref{Th6}, for $\sigma_{n+2}=\pm$ we have
\begin{equation}
{\cal T}_{n+1}^{j_n+1}=\Big\{O\cap\Big(\bigcup _{{\delta_0}...{\delta_{n+1}}} N_{\delta_0...\delta_{n+1}}^{j_{n+1}}\Big) \quad /\quad O\in {\cal T}_{n+2}^{j_{n+2}}\Big\},
\end{equation}
then we have to find $O\in {\cal T}_{n+2}^{j_{n+2}}$ such that $\Omega=O\cap\Big(\bigcup _{{\delta_0}...{\delta_{n+1}}} N_{\delta_0...\delta_{n+1}}^{j_{n+1}}\Big)$.

In this purpose, let us consider the set
\begin{equation}\label{O3}
O=\bigcup _{{\delta_0}...{\delta_{n+2}}} O_{\delta_0...\delta_{n+2}}^{\sigma_0\ldots\sigma_{n+2}}
\end{equation}
defined by

\begin{equation}\label{O4}
O_{\delta_0...\delta_{n+2}}^{\sigma_0...\sigma_{n+2}}=
\left\{
  \begin{array}{ll}
    \Omega_{\delta_0...\delta_{n}}^{\sigma_0...\sigma_{n}}& \forall \delta_0,\ldots,\delta_n\quad  \hbox{and for} \quad \delta_{n+1}=0, \ \delta_{n+2}=0\\
    \emptyset & \forall \delta_0,\ldots,\delta_n\quad  \hbox{and for} \quad \delta_{n+1}\neq0, \ \delta_{n+2}=0\\
    \emptyset & \forall \delta_0,\ldots,\delta_{n+1}\quad  \hbox{and for} \quad \delta_{n+2}\neq0\\
  \end{array}
\right.
\end{equation}

$O\in {\cal T}_{n+2}^{j_{n+2}}$, indeed:
\begin{itemize}
  \item {for all $\delta_0,\ldots,\delta_{n+1}$ and for $\delta_{n+2}\neq0$,\quad $\emptyset$ is an open set of
$\bigcup _{{\delta_0}...{\delta_{n+1}}{\delta_{n+2}}} N_{\delta_0...\delta_{n+1}\delta_{n+2}}^{j_{n+2}}$.}
  \item {for all $\delta_0,\ldots,\delta_n$, for $\delta_{n+1}\neq0$, and $\delta_{n+2}=0$,\ $\emptyset$ is an open set of
$\bigcup _{{\delta_0}...{\delta_{n+1}}{\delta_{n+2}}} N_{\delta_0...\delta_{n+1}0}^{j_{n+2}}$. }
  \item {for all $\delta_0,\ldots,\delta_n$, $\delta_{n+1}=0$ and $ \delta_{n+2}=0$, $\Omega_{\delta_0...\delta_{n}}^{\sigma_0...\sigma_{n}}$ is an open set of
$N_{\delta_0...\delta_{n}}^{j_n}$. By applying the property \ref{pr1} twice, we have the following identification
$$N_{\delta_0...\delta_{n}}^{j_n}=N_{\delta_0...\delta_{n}0}^{j_{n+1}}=N_{\delta_0...\delta_{n}00}^{j_{n+2}},$$
then $\Omega_{\delta_0...\delta_{n}}^{\sigma_0...\sigma_{n}}$ is an open set of
$N_{\delta_0...\delta_{n}00}^{j_{n+2}}$.}
\end{itemize}

Let us prove now that $\Omega=O\cap\Big(\bigcup _{{\delta_0}...{\delta_{n+1}}} N_{\delta_0...\delta_{n+1}}^{j_{n+1}}\Big)$:
by (\ref{O3}) and (\ref{O4}),
$$O\cap\Big(\bigcup _{{\delta_0}...{\delta_{n+1}}} N_{\delta_0...\delta_{n+1}}^{j_{n+1}}\Big)=
\Big(\bigcup _{{\delta_0}...{\delta_{n+2}}} O_{\delta_0...\delta_{n+2}}^{\sigma_0\ldots\sigma_{n+2}}\Big)
\cap\Big(\bigcup _{{\delta_0}...{\delta_{n+1}}} N_{\delta_0...\delta_{n+1}}^{j_{n+1}}\Big)$$

$$=\Big((\bigcup_{\delta_{n+2}\neq0}\bigcup_{{\delta_0}...{\delta_{n+1}}}\emptyset)
\cup(\bigcup_{{{\delta_{n+1}}\atop {\delta_{n+2}}}{\neq\atop =}{0\atop 0}
}\bigcup_{\delta_0...\delta_n}\emptyset)
\cup(\bigcup_{{{\delta_{n+1}}\atop {\delta_{n+2}}}{=\atop =}{0\atop 0}}\bigcup_{\delta_0...\delta_n}\Omega_{\delta_0...\delta_{n}}^{\sigma_0...\sigma_{n}})\Big)
\cap\Big(\bigcup _{{\delta_0}...{\delta_{n+1}}} N_{\delta_0...\delta_{n+1}}^{j_{n+1}}\Big)$$

$$=\Big(\bigcup_{{{\delta_{n+1}}\atop {\delta_{n+2}}}{=\atop =}{0\atop 0}}\bigcup_{\delta_0...\delta_n}\Omega_{\delta_0...\delta_{n}}^{\sigma_0...\sigma_{n}}\Big)
\cap\Big(\bigcup _{{\delta_0}...{\delta_{n+1}}} N_{\delta_0...\delta_{n+1}}^{j_{n+1}}\Big)$$
Since
$$\Big(\bigcup_{{{\delta_{n+1}}\atop {\delta_{n+2}}}{=\atop =}{0\atop 0}}\bigcup_{\delta_0...\delta_n}\Omega_{\delta_0...\delta_{n}}^{\sigma_0...\sigma_{n}}\Big)
\cap\Big(\bigcup _{{\delta_0}...{\delta_{n+1}}} N_{\delta_0...\delta_{n+1}}^{j_{n+1}}\Big)=\emptyset\quad \hbox{for } \delta_{n+1}\neq0,$$
then we have

$$O\cap\Big(\bigcup _{{\delta_0}...{\delta_{n+1}}} N_{\delta_0...\delta_{n+1}}^{j_{n+1}}\Big)=
\Big(\bigcup_{\delta_0...\delta_n}\Omega_{\delta_0...\delta_{n}}^{\sigma_0...\sigma_{n}}\Big)
\cap\Big(\bigcup _{{\delta_0}...{\delta_{n}}0} N_{\delta_0...\delta_{n}0}^{j_{n+1}}\Big)$$

$$=\Omega\cap\Big(\bigcup _{{\delta_0}...{\delta_{n}}0} N_{\delta_0...\delta_{n}0}^{j_{n+1}}\Big).$$
By lemma \ref{L1},

$$O\cap\Big(\bigcup _{{\delta_0}...{\delta_{n+1}}} N_{\delta_0...\delta_{n+1}}^{j_{n+1}}\Big)
=\Omega\cap\Big(\bigcup _{{\delta_0}...{\delta_{n}}} N_{\delta_0...\delta_{n}}^{j_{n}}\Big)=\Omega,$$
that proves the result.
Therefore, the family (\ref{Ft1}) is a fractal family of topological spaces.

\rightline\Box

We have the following immediate consequences:

\begin{cor}\label{Co1}
For all $n\geq0$ and for $j_n=\sigma_0...\sigma_{n}\in \Lambda_n$,
the topology ${\cal T}_{n}^{j_{n}}$ is weaker than the topology ${\cal T}_{n+1}^{j_{n+1}}$ for $j_{n+1}=\sigma_0...\sigma_{n}\sigma_{n+1}\in \Lambda_{n+1}$.
\end{cor}
\pesp

\begin{cor}
The weakest topology of the fractal family of topological spaces
\begin{equation}
\Big(\bigcup _{{\delta_0}...{\delta_{n}}} N_{\delta_0...\delta_{n}}^{j_{n}},\ {\cal T}_{n}^{j_{n}}\Big)_{{j_n\atop n}{\in\Lambda_n\atop\geq0}}
\end{equation}
is the topology ${\cal T}_{0}^{\sigma_0}$ for $\sigma_0=\pm$.
\end{cor}

\begin{thm}\label{Th2}
The fractal manifold  ${\cal M}$ is locally homeomorphic via a family of local homeomorphisms to the fractal topological space
$\Big(\bigcup _{{\delta_0}...{\delta_{n}}} N_{\delta_0...\delta_{n}}^{j_{n}}\Big)_{{j_n\atop n}{\in\Lambda_n\atop\geq0}}$.
\end{thm}

\ni{\bf Proof.}  By  Theorem \ref{Th1}, the fractal manifold $\cal M$ is locally homeomorphic via a family of local homeomorphisms to the family of diagonal
topological spaces given by (\ref{Ft1}). The use of theorem \ref{Th5} and definition \ref{Def4} ends the proof.

\rightline\Box

\begin{prop}\label{prop8}
Let ${\cal M}$ be a fractal manifold locally homeomorphic via a family of homeomorphisms to the fractal topological space
$\Big(\bigcup _{{\delta_0}...{\delta_{n}}} N_{\delta_0...\delta_{n}}^{j_n}\Big)_
{{j_n\atop n}{\in\Lambda_n\atop\geq0}}$
endowed with the fractal topology
$\Big({\cal T}_{n}^{j_n}\Big)_
{{j_n\atop n}{\in\Lambda_n\atop\geq0}}$, then a local chart of ${\cal M}$ at the $step(n)$ is given by a $(Card\ \Lambda_{n+1} +1)$-uplet.
\end{prop}

\ni{\bf Proof.} For $n=0$, using theorem $\ref{Th1}$ at the $step(0)$, we have one diagram with a double homeomorphism,
and the local chart is given by a triplet $\ (\Omega, \varphi_1, T_1\circ \varphi_1)$. For $n=0$, $\Lambda_0=\{\sigma_0\ /\ \sigma_0=+,\ \sigma_0=-\}$, then
$Card\ \Lambda_0=2=2^1$. Then the local chart at the $step(0)$ is given by a $(2^1+1)$-uplet, that is to say a $(Card\ \Lambda_0+1)$-uplet.

At the $step(1)$, using theorem \ref{Th1}, each local homeomorphism at the $step(0)$ becomes a double local homeomorphism, and the local chart is given by a $(2^2+1)$-uplet
$$\ (\Omega, \varphi_2\circ\varphi_1, T_2\circ \varphi_2\circ\varphi_1,\varphi_3\circ T_1\circ \varphi_1, T_3\circ\varphi_3\circ T_1\circ \varphi_1).$$
For $n=1$, $\Lambda_1=\{\sigma_0\sigma_1\ /\ \sigma_0=\pm,\ \sigma_1=\pm\},$
then $Card\ \Lambda_1=2\times Card\ \Lambda_0=2^2$.
Then the local chart is given by a $(2^2+1)$-uplet, that is to say a $(Card\ \Lambda_1+1)$-uplet.

By induction over $n\geq0$, suppose that the property is true at the $step(n)$, that is to say: a local chart of $\cal M$ is given by a $(Card\ \Lambda_n+1)$-uplet,
 where $Card\ \Lambda_n=2^{n+1}$.

Let us prove the property at the $step(n+1)$: using theorem \ref{Th1}, each local homeomorphism at the $step(n)$
becomes a double local homeomorphism at the $step(n+1)$, and the local chart is given by a $(2^{n+1}\times 2 +1)$-uplet, that is to say a $(2^{n+2} +1)$-uplet.

Since $\Lambda_{n+1}=\{\sigma_0\ldots\sigma_n\sigma_{n+1}\ /\ \sigma_0=\pm,\ \ldots\ ,\sigma_n=\pm,\sigma_{n+1}=\pm\}$ and $Card\ \Lambda_n=2^{n+1}$, then
$Card\ \Lambda_{n+1}=2Card\ \Lambda_{n}=2^{n+2}$.

Then a local chart at the $step(n+1)$ is given by a $(Card\ \Lambda_{n+1}+1)$-uplet.

\rightline\Box

\section{Conclusion}

We can wonder in which domain a fractal topology could be applied. The idea
of a space-time that fluctuates at the microscopic level appeared in the 1960s \cite{W}.
Moreover it is known that the general relativity is unable to describe the topology of the universe.
Accordingly  new perspectives could be open by modeling the space-time by a fractal manifold together with its fractal topology. It could offer
several valuable insights into a concrete description of an expanding space-time and enrich the field of the cosmic topology.

% Create the reference section using BibTeX:
%\begin{acknowledgements}
%If you'd like to thank anyone, place your comments here
%and remove the percent signs.
%\end{acknowledgements}

% BibTeX users please use one of
%\bibliographystyle{spbasic}      % basic style, author-year citations
\bibliographystyle{spmpsci}      % mathematics and physical sciences
%\bibliographystyle{spphys}       % APS-like style for physics
%\bibliography{}   % name your BibTeX data base

\end{document}